\documentclass[a4paper,11pt,leqno]{amsart}
\usepackage{hyperref}
\usepackage{amsmath}
\usepackage{amssymb}
\usepackage{amsthm}
\usepackage{amsopn}
\usepackage{amscd}
\usepackage{textcomp}
\usepackage[all]{xy}
\usepackage{lscape}
\usepackage{verbatim}

\newcommand{\frS}{\mathfrak{S}}

\newcommand{\frU}{\mathfrak{U}}

\newcommand{\frZ}{\mathfrak{Z}}

\newcommand{\frg}{\mathfrak{g}}
\newcommand{\frh}{\mathfrak{h}}

\newcommand{\frk}{\mathfrak{k}}
\newcommand{\frl}{\mathfrak{l}}

\newcommand{\frqqq}{\mathfrak{q}}

\newcommand{\frt}{\mathfrak{t}}
\newcommand{\fru}{\mathfrak{u}}

\newcommand{\bbC}{\mathbb{C}}

\newcommand{\bbH}{\mathbb{H}}

\newcommand{\bbN}{\mathbb{N}}

\newcommand{\bbR}{\mathbb{R}}

\newcommand{\bbZ}{\mathbb{Z}}

\newcommand{\caB}{\mathcal{B}}

\newcommand{\caD}{\mathcal{D}}

\newcommand{\caI}{\mathcal{I}}

\newcommand{\caM}{\mathcal{M}}

\newcommand{\caP}{\mathcal{P}}

\newcommand{\caR}{\mathcal{R}}

\newcommand{\caU}{\mathcal{U}}

\newcommand{\caX}{\mathcal{X}}

\def\n{{\mathbb N}}

\newcommand{\sgn}{\mathrm{sgn}}

\theoremstyle{plain}
\newtheorem{thm}{Theorem}[section]
\newtheorem{lemme}[thm]{Lemma}
\newtheorem{cor}[thm]{Corollary}
\newtheorem{prop}[thm]{Proposition}

\theoremstyle{definition}

 \newtheorem{defi}[thm]{Definition}
 
\newtheorem{rmq}[thm]{Remark}

\newcommand{\Hom}{\mathrm{Hom}}
\newcommand{\im}{\mathrm{Im}\, }

\newcommand{\supp}{\mathbf{Supp }\,  }
\newcommand{\tr}{\mathrm{Tr}\, }

\newcommand{\C}{\mathbf{C}}
\newcommand{\GL}{\mathbf{GL}}
\newcommand{\Irr}{\mathbf{Irr}}
\newcommand{\Lg}{\mathbf{Lg}}
\newcommand{\SU}{\mathbf{SU}}
\newcommand{\Or}{\mathbf{O}}
\newcommand{\U}{\mathbf{U}}
\newcommand{\Sp}{\mathbf{Sp}}

\newcommand{\lra}{\leftrightarrow}

\newcommand{\JL}{\mathbf{LJ}}
\newcommand{\ccc}{\mathbf{C}}

\newcommand{\LJ}{\mathbf{LJ}}

\newcommand{\fff}{\mathcalF}

\def \pr {\noindent \underline{\sl Proof}. }

\begin{document}


\title[Global Jacquet-Langlands correspondence]
{Unitary dual of $\GL(n)$ at archimedean places and global Jacquet-Langlands correspondence}
\author{I. A Badulescu  and  D. Renard}
\address{
Centre de Mathématiques 
Laurent Schwartz
\\ Ecole Polytechnique }
\date{\today}

\begin{abstract}
In \cite{Ba1}, results about the global Jacquet-Langlands correspondence, (weak and strong) multiplicity-one theorems
and the classification of automorphic representations for inner forms of the  general linear group over a
 number field    are established, under the  condition that the local inner forms are split at  archimedean places.
In this paper, we extend the main local results of \cite{Ba1} to archimedean places so that this assumption 
can be removed. Along the way, we collect several results about the unitary dual of general linear groups over 
$\bbR$, $\bbC$ or $\bbH$ of independent interest.
\end{abstract}

\maketitle

\tableofcontents


\section{Introduction}

In \cite{Ba1}, results about the global Jacquet-Langlands correspondence, (weak and strong)   
multiplicity-one theorems and the classification of automorphic representations for inner forms of the  
general linear group over a number field    are established, under the  condition that the local inner forms
 are split at  archimedean places. The main goal of
 this paper is to remove this hypothesis. The paper consists in two parts. In the first part, we extend 
 the main local results of \cite{Ba1} to archimedean places. In the second part, we explain how to use
 these local results to establish the global results in their full generality. Along the way, we collect 
several results about the unitary  dual of general linear groups over $\bbR$, $\bbC$ or $\bbH$ of independent 
interest. Let us now explain in more details the content of this paper.

\subsection{Some notation}
Let $A$ be one of  the division algebras $\bbR$, $\bbC$ or $\bbH$. If $A=\bbR$ or $A=\bbC$ and $n\in \bbN^\times$,
 we denote by $\det$ the determinant map on $\GL(n,A)$ (taking values in $A$). If $A=\bbH$, let $RN$ be 
the reduced norm map on $\GL(n,\bbH)$ (taking values in $\bbR^\times_+$).

If $n \in \bbN$ and   $\sum_{i=1}^s n_i=n$ is a partition of  $n$, the group  
$\GL({n_1},A)\times \GL({n_2},A)\times ...\times \GL({n_s},A)$
is identified with the subgroup of  $\GL(n,A)$ of bloc diagonal matrices of size 
 $n_1,\ldots, n_s$. Let us denote  $G_{(n_1,\ldots, n_s)}$ this subgroup and 
$P_{(n_1,\ldots, n_s)}$ the  parabolic subgroup of  $GL(n,A)$ containing 
  $G_{(n_1,\ldots, n_s)}$ and the  Borel subgroup of invertible upper triangular matrices.
For $1\leq i\leq s$ let  $\pi_i$ be an admissible representation of $\GL({n_i},A)$ of finite length. 
We write then $\pi_1\times \pi_2\times ...\times \pi_s$ for the representation induced from the representation
$\pi_1\otimes \pi_2\otimes ...\otimes \pi_s$ of $G_{(n_1,\ldots, n_s)}$ with respect to $P_{(n_1,\ldots, n_s)}$. 
We will also use this notation  for the image of representations in the Grothendieck group of virtual characters, 
 which makes the above product commutative.  We also  often do not distinguish between  a representation
and  its  isomorphy class  and write  ``equal'' for ``isomorphic''.

\subsection {Classification of unitary representations}

We recall first Tadi\'c  classification of the unitary dual of the groups $\GL(n,\bbR)$ and  $\GL(n,\bbC)$,
 following \cite{Ta5}. The classification is similar to the one for non archimedean local fields
(\cite{Ta1},\cite{Ta2}) and is explained in details in Section \ref{UD}. Part
 of the arguments do not appear in the literature in the case of  $\GL(n,\bbH)$, so we give the complete 
proofs in Section \ref{U3}, \ref{U1ARC} and \ref{VOG}, using Vogan's classification \cite{Vog2}.

Let $X_\bbC$ be the set of unitary characters of $\bbC^\times$. If $\chi\in X_\bbC$, $n\in \bbN^\times$ 
let $\chi_n$ be the character $\chi\circ\det$ of $\GL(n,\bbC)$. Let $\nu_n$ be the character of
 $\GL(n,\bbC)$ given by the square of the module of the determinant. If $\sigma$ is a 
representation of $\GL(n,\bbC)$ and $\alpha\in \bbR$, write $\pi(\sigma,\alpha)$ for
 the representation $\nu_n^\alpha\sigma\times \nu_n^{-\alpha}\sigma$ of
$\GL(2n,\bbC)$.
Set 
$$\caU_\bbC =\{\chi_n,\pi(\chi_n,\alpha)\ |\ \ \chi\in X_\bbC,\, n\in \bbN^\times,\,
 \alpha\in ]0,\frac{1}{2}[\}.$$

Let $X_\bbR$ be the set of unitary characters of $\bbR^\times$. Let  sgn denote the sign  character.
If $\chi\in X_\bbR$, $n\in \bbN^\times$ let $\chi_n$ be the character $\chi\circ\det$ of $\GL(n,\bbR)$ and
 $\chi'_n$ the character $\chi\circ RN$ of $\GL(n,\bbH)$. For fixed $n$,
the map $\chi\mapsto \chi_n$ is an isomorphism from the group of unitary characters of
 $\bbR^\times$ to the group of unitary characters of $\GL(n,\bbR)$, while  $\chi\mapsto \chi'_n$ 
is a surjective map from the group of unitary characters of $\bbR^\times$ to the group of unitary characters
 of $\GL(n,\bbH)$, with kernel $\{1,\sgn\}$.

Let $\nu_n$ (resp. $\nu'_n$) be the character of $\GL(n,\bbR)$ (resp. $\GL(n,\bbH)$) given
 by the absolute  value (resp. the reduced norm) of the determinant.
If $\sigma$ is a representation of $\GL(n,\bbR)$ (resp. $\GL(n,\bbH)$) and $\alpha\in \bbR$, 
write $\pi(\sigma,\alpha)$ for the representation $\nu_n^\alpha\sigma\times \nu_n^{-\alpha}\sigma$ of
$\GL(2n,\bbR)$ (resp. the representation ${\nu'}_n^\alpha\sigma\times {\nu'}_n^{-\alpha}\sigma$ of
$\GL(2n,\bbH)$).

Let $D^u_2$ be the set of isomorphy classes of square integrable (modulo center) representations of $\GL(2,\bbR)$.
 For $\delta\in D^u_2$ and $k\in \bbN^\times$, write $u(\delta,k)$ for the Langlands quotient of the representation
 $\nu_2^{\frac{k-1}{2}}\delta\times\nu_2^{\frac{k-3}{2}}\delta\times\nu_2^{\frac{k-5}{2}}\delta\times...
\times\nu_2^{-\frac{k-1}{2}}\delta$. Then $u(\delta,k)$ is a representation of $\GL({2k},\bbR)$. 
Set
\begin{align*}
& \caU_\bbR =\{\chi_n, \pi(\chi_n,\alpha)\ |\ \ \chi\in X_\bbR,\, n\in \bbN^\times,\, 
\alpha\in ]0,\frac{1}{2}[\}\\
& \cup  \{u(\delta,k),\pi(u(\delta,k),\alpha)\ |\ \ \delta\in D^u_2,\, k\in \bbN^\times,\ \alpha\in ]0,\frac{1}{2}[\}.
\end{align*}


Let now $D$ be the set of unitary representations of $\bbH^\times$ {\it which are not one-dimensional}.
For $\delta\in D$ and $k\in \bbN^\times$, write $u(\delta,k)$ for the Langlands quotient of the 
representation ${\nu'}_1^{\frac{k-1}{2}}\delta\times{\nu'}_1^{\frac{k-3}{2}}\delta\times{\nu'}_1^{\frac{k-5}{2}}
\delta\times...\times{\nu'}_1^{-\frac{k-1}{2}}\delta$. Then $u(\delta,k)$ is a representation of $\GL(k,\bbH)$. 

Set
\begin{align*}
&\caU_\bbH =\{\chi'_n,\pi(\chi'_n,\alpha)\ |\ \ \chi\in X_\bbR,\, n\in \bbN^\times,\, \alpha\in ]0,1[\}\\
&\cup\ \ 
\{u(\delta,k),\pi(u(\delta,k),\alpha)\ |\ \ \delta\in D,\, k\in \bbN^\times\, \alpha\in ]0,\frac{1}{2}[\}.
\end{align*}

\begin{thm}
For $A=\bbC,\bbR,\bbH$,  any representation in $\caU_A$ is irreducible and unitary, any product of representations in
 $\caU_A$ is irreducible and unitary, and any irreducible unitary representation $\pi$ of 
$\GL(n,A)$ can be written as a product of elements in $\caU_A$. Moreover, $\pi$ determines the factors of the product 
(up to permutation).  
\end{thm}
 Notice the two different ranges for the  the possible values of $\alpha$ in the case   $A=\bbH$.

\subsection{Jacquet-Langlands for unitary representations}
Any element in $\GL(n,\bbH)$ has a characteristic polynomial of degree $2n$ with coefficients  in $\bbR$.
We say that two elements $g\in \GL(2n,\bbR)$ and $g'\in \GL(n,\bbH)$ correspond (to each other) 
if they have the same characteristic polynomial and this polynomial has distinct roots in 
$\bbC$ (this last condition means that $g$ and $g'$ are regular semisimple). We then write $g\lra g'$.

Let $\C$ denote the Jacquet-Langlands  correspondence between irreducible square integrable 
representations of $\GL(2,\bbR)$ and irreducible unitary representations of $\bbH^\times$ (\cite{JL}). 
This correspondence can be extended to a correspondence $|\LJ|$ between {\it all} irreducible
 unitary representations of $\GL(2n,\bbR)$ and $\GL(n,\bbH)$ (it comes from a ring morphism 
$\LJ$ between the respective Grothendieck groups, defined in Section \ref{JL1}, which explain the notation).
 In what follows, it is understood that each time we write the relation $|\LJ|(\pi)=\pi'$ for $\pi$ and $\pi'$ 
representations of   $\GL(2n,\bbR)$ and $\GL(n,\bbH)$ respectively, then  $\pi$ and $\pi'$ satisfy the 
character relation $\Theta_\pi(g)=\varepsilon(\pi) \Theta_\pi'(g')$ for all $g\lra g'$ 
where $\varepsilon(\pi)$ is an explicit  sign ($\pi$ clearly determines $\pi'$ and $\varepsilon$). 
The correspondence $|\LJ|$ for unitary representations is  given first   on elements 
in $\caU_\bbR$:\\

(a) $|\LJ|(\chi_{2n})=\chi'_n$ and $|\LJ|(\pi(\chi_{2n},\alpha))=\pi(\chi'_n,\alpha)$ 
for all $\chi\in X_\bbR$ and $\alpha\in ]0,\frac{1}{2}[$;\\

(b) If $\delta\in D^u_2$ is such that $\C(\delta)$ is in $D$ (i.e. is not one-dimensional) 
then $|\LJ|(u(\delta,k))=|\LJ|(\C(\delta),k)$ and $\LJ(u(\delta,k),\alpha)= \pi(u(\C(\delta),k),\alpha)$ 
for all $\alpha\in ]0,\frac{1}{2}[$;\\

(c) If $\delta\in D^u_2$ is such that $\C(\delta)$ is a one-dimensional representation $\chi'_1$, then

--- $|\LJ|(u(\delta,k))=\pi(\chi'_\frac{k}{2},\frac{1}{2})$, $|\JL|(\pi(u(\delta,k),\alpha))=
\pi(\pi(\chi'_\frac{k}{2},\frac{1}{2}),\alpha)$ if $k$ is even and $\alpha\in ]0,\frac{1}{2}[$. \\

--- $|\JL|(u(\delta,k))=\chi'_\frac{k+1}{2}\times\chi'_\frac{k-1}{2}$, 
 $|\JL|(\pi(u(\delta,k),\alpha))=\pi(\chi'_\frac{k+1}{2},\alpha)\times \pi(\chi'_\frac{k-1}{2},
 \alpha)$  if $k\neq 1$  odd and $\alpha\in ]0,\frac{1}{2}[$. \\

--- $|\JL|(\delta)=\chi'_1$,  $|\JL|(\pi(\delta,\alpha))=\pi(\chi'_1,\alpha)$,  $\alpha\in ]0,\frac{1}{2}[$,.\\

Let $\pi$ be an irreducible unitary representation of $\GL(2n,\bbR)$. If writing $\pi$ 
as a product of elements in  $\caU_\bbR$ involves a factor  not listed in (a), (b) or (c) 
it is easy to show that $\pi$ has a character which vanishes on elements which correspond to 
elements of $\GL(n,\bbH)$, and we set $|\JL|(\pi)=0$. If all the factors $\sigma_i$  of $\pi$ are in 
 (a) (b) (c) above,   $|\JL|(\pi)$ is the product of the  $ |\JL|(\sigma_i)$ (an irreducible unitary
 representation of  $\GL(n,\bbH)$). 
Elements of $\caU_\bbR$ not listed at (a) (b) (c) are of type $\chi$ or 
$\pi(\chi,\alpha)$, with $\chi$ a character of some $\GL(k,\bbR)$ and $k$ odd.

Notice that some unitary irreducible representations of $\GL(n,\bbH)$ are not in the image 
of this map (if $n\geq 2$). For instance, when $\chi\in X_\bbR$ and  $\frac{1}{2}<\alpha <1$, 
then both $\pi(\chi_2,\alpha)$ and $\pi(\chi'_1,\alpha)$ are irreducible and correspond 
to each other by the character relation, but $\pi(\chi'_1,\alpha)$ is unitary while
 $\pi(\chi_2,\alpha)$ is not. Using the classification of unitary representations for
 $\GL(4,\bbR)$ and basic information from the infinitesimal character, it is clear 
that no (possibly other) unitary representation of $\GL(4,\bbR)$ has matching 
character with $\pi(\chi'_1,\alpha)$.

As a consequence of the above results, we get:

\begin{thm}\label{unit}
Let $u$ be a unitary irreducible  representation of $\GL(2n,\bbR)$. Then either 
the character $\Theta_u$ of $u$ vanishes on the set of elements of $\GL(2n,\bbR)$ which
 correspond to some element of $\GL(n,\bbH)$, or there exists
a unique irreducible unitary (smooth) representation $u'$ of $\GL(n,\bbH)$ such that
$$\Theta_u(g)=\varepsilon(u) \Theta_{u'}(g')$$
for all $g\lra g'$, where $\varepsilon(u)\in \{-1,1\}$.
\end{thm}

The above results are proved in Section \ref{JaLa} and are based on the fact that  $\GL(2n,\bbR)$
and  $\GL(n,\bbH)$ share Levi subgroups (of $\theta$-stable parabolic subgroups, the ones used in cohomological
 induction (\cite{KV})) which are products of  $\GL(n_i,\bbC)$. The underlying principle (a nice instance of Langlands'
functoriality) is that the Jacquet-Langlands morphism  $\JL$ commutes with cohomological induction. The same 
principle, with Kazhdan-Patterson lifting instead of Jacquet-Langlands correspondence, was already used in \cite{AH}.

\subsection{Character identities and ends of complementary series}

In Section \ref{ends}, we give the composition series of ends of complementary series in most cases. 
This is not directly related to the main purpose of the paper, which is the global theory 
of the second part, but it solves some old conjectures of Tadi\'c which are important in understanding the 
topology of the unitary dual of the groups $\GL(n,A)$, $A=\bbR,\bbC,\bbH$. The starting point is 
Zuckerman formula for the trivial representation of $\GL(n,A)$. Together with cohomological induction, it gives
character formulas for unitary representations of the groups $\GL(n,A)$. In case $A=\bbC$, Zuckerman 
formula is given by a determinant (see formula (\ref{Zu2})), and the Lewis Carroll identity of \cite{CR}
allows us to deduce  formulas (\ref{boutC}), (\ref{boutR}),  (\ref{boutR1}), (\ref{boutH}), (\ref{boutH1}) 
for the ends of complementary series.

\subsection{Global results}
\def\aa{{\mathbb A}}
\def\bc{\backslash}
\def\lg{L^2(G(F)Z(\aa)\bc G(\aa))}
\def\o{\omega}
\def\lgg{L^2(G'(F)Z(\aa)\bc G'(\aa);\o)}
\def\lra{\leftrightarrow}
\def\ccc{{\bf C}}
\def\jlr{{\bf JL}_r}
\def\ljr{{\bf LJ}_r}

\def\s{{\mathfrak S}}
\numberwithin{equation}{section}
\leftskip -1cm
\rightskip -1cm
\def\a{\alpha}
\def\b{{\mathcal B}}
\def\d{{\mathcal D}}
\def\e{\varepsilon}
\def\f{{\mathcal F}}
\def\i{{\bf i}}
\def\k{\{1,2,...,k\}}
\def\l{{\bf l}}
\def\n{\mathbb N}
\def\r{{\mathbb R}}
\def\s{\sigma}
\def\z{\mathbb Z}

\def\cc{{\mathbb C}}
\def\ccc{{\bf C}}
\def\lra{\leftrightarrow}
\def\ski{{\sum_{i=1}^k}}
\def\ki{{_{i=1}^k}}
\def\rrr{{\mathcal R}}

\def\lj{{\bf LJ}}
\def\sgls{standard Levi subgroup}
\def\sglss{standard Levi subgroups}
\def\bc{\backslash}
\def\cusp{{\mathcal C}}
\def\tr{{\rm tr}}

\def\g{{\bf G}}

\def\fff{{\mathcal F}}
Let ${\mathcal F}$ be a global field of characteristic zero and ${\mathcal D}$ a
central division algebra over ${\mathcal F}$ of
dimension $d^2$.
Let $n\in\n^*$. Set $A'=M_{n}({\mathcal D})$. For each place $v$ of ${\mathcal F}$ let
${\mathcal F}_v$ be the completion of ${\mathcal F}$ at $v$ and set $A'_v=A'\otimes {\mathcal F}_v$.
For every place $v$ of ${\mathcal F}$, $A'_v$ is isomorphic to the matrix algebra 
$M_{r_v}({\mathcal D}_v)$ for some positive number $r_v$ and some central division
algebra ${\mathcal D}_v$ of dimension
$d_v^2$ over ${\mathcal F}_v$ such that $r_v d_v=nd$. We will fix once and for all an isomorphism and
identify these two algebras. Let $V$ be the (finite) set of places
where $M_n({\mathcal D})$ is not split (i.e. $d_v\neq 1$).

Let $G'({\mathcal F})$ be the group $A'^\times=\GL(n,{\mathcal D})$.
For every place $v\in V$, set $G'_v= A'^\times_v= \GL(r_v, {\mathcal D}_v)$ and 
$G_v=\GL(n, \fff_v)$. For a given place $v$ (clear from the context)
 write $g\lra g'$ if $g\in G_v$ and $g'\in G_v'$ are regular semisimple and have equal characteristic polynomial.

If $v\notin V$, the groups $G_v$ and $G'_v$ are isomorphic and we  fix  once and for all an isomorphism
which allows us to identify them (as every  automorphism of $G_v$ or $G'_v$ is inner, the 
choice of this isomorphism will not be relevant when working with equivalence classes of representations).

Theorem \ref{unit} has been proved in the $p$-adic case too (\cite{Ta4}, \cite{Ba1}).
 So, if $v\in V$, with the same notation and conventions in the $p$-adic case as in
 the archimedean case :  

\begin{thm}\label{unit2}
Let $u$ be a unitary irreducible smooth representation of $G_v$. Then we have one and
 only one of the two following :

{\rm (i)} the character $\Theta_u$ of  $u$ vanishes on the set of elements of $G_v$ which
  correspond to elements  of $G'_v$,

{\rm (ii)} there exists
a unique unitary smooth irreducible representation $u'$ of $G'_v$ such that
$$\Theta_u(g)=\varepsilon(u) \Theta_{u'}(g')$$
for any $g\lra g'$, where $\varepsilon(u)\in \{-1,1\}$.
\end{thm}

In the second case (ii) we say $u$ is {\bf compatible}. We denote the map $u\mapsto u'$
 defined on the set of compatible (unitary) representations by $|\LJ_v|$.

Let $\aa$ be the ring of adeles of ${\mathcal F}$, $G'(\aa)$.
The group  $G'(\fff)$ (resp. $G(\fff)$) is a discrete  subgroup of $G'(\aa)$  (resp. $G(\aa)$).
The centers of $G'$ and $G$ consist of scalar nonzero  matrices  and can be identified, so both will be 
denoted by $Z$. 

We endow these local and global groups with measures as in \cite{AC}
and for every unitary smooth character $\o$ of $Z(\aa)$ trivial on
$Z({\fff})$, we let $L^2(G'(\fff)Z(\aa)\bc G'(\aa);\o)$
be the space of functions $f$ defined on $G'(\aa)$ with values in
$\cc$ such that

i) $f$ is left invariant under $G'({\mathcal F})$,

ii)  $f(zg)=\o(z)f(g)$ for all $z\in Z(\aa)$ and all $g\in G'(\aa)$,

iii) $|f|^2$ is integrable over $G'({\mathcal F})Z(\aa)\bc G'(\aa)$.\\
Let us denote by  $R'_\o$ the representation of $G'(\aa)$ on $\lgg$ by
right translations. A {\it discrete series of $G'(\aa)$} is the equivalence class of  an
irreducible subrepresentation of $R'_\o$  for some smooth unitary character $\o$ of $Z(\aa)$ trivial on
$Z({\mathcal F})$. Then $\o$ is the central character of $\pi$. Let $R'_{\o,disc}$ be the 
 subrepresentation of $R'_\o$ generated by irreducible subrepresentations. 
It is known that a discrete series representation of  $G'(\aa)$  appears with finite 
multiplicity  in $R'_{\o,disc}$.

Similar definitions and statements can be made  with $G$ instead of $G'$, with obvious notation. 
Every discrete series $\pi$ of $G'(\aa)$ (resp. $G(\aa)$) is ``isomorphic'' to
a restricted Hilbert  tensor product of irreducible unitary smooth
representations $\pi_v$ of the groups $G'_v$ (resp. $G_v$) - see  \cite{Fl} for a precise statement and proof. 
The local components $\pi_v$ are determined by $\pi$.

Denote $DS$ (resp. $DS'$) the set of discrete series of $G(\aa)$ (resp. $G'(\aa)$).
 Let us  say that a discrete series $\pi$ of $G(\aa)$ is
${\mathcal D}${\bf -compatible} if $\pi_v$ is compatible for all places $v\in V$.

\begin{thm}\label{corr}
{\rm (a)} There exists a unique map $\g:DS'\to DS$ such that for every $\pi'\in DS'$, if $\pi=\g(\pi')$, one has

--- $\pi$ is $\caD$-compatible,

--- if $v\notin V$, then $\pi_v=\pi'_v$ and

--- if $v\in V$, then $|\JL_v|(\pi_v)=\pi'_v$.\\
The map $\g$ is injective. The image of $\g$ is the set
of all ${\mathcal D}$-compatible discrete series of $G(\aa)$.

{\rm (b)} If $\pi'\in DS'$, then
the multiplicity of $\pi'$ in the discrete spectrum is one (Multiplicity One Theorem).

{\rm (c)} If $\pi', \pi''\in DS'$, if $\pi'_v\simeq\pi''_v$ for almost all $v$, 
then $\pi'=\pi''$ (Strong Multiplicity One Theorem).
\end{thm}

With $\mathcal D$ fixed, we need now to consider all possible  $n\in \bbN^*$ at the same time and we add a  subscript 
in the notation :   $A'_n=M_n(\fff)$, $A_n=M_n({\mathcal D})$,  $G_n$, $G'_n$, $DS_n$, $DS'_n$...
We recall  the Moeglin-Waldspurger classification of the residual
spectrum for the groups $G_n(\aa)$, $n\in \n^*$. Let $\nu$ be the
character of $G_n(\aa)$ or $G'_n(\aa)$ given by the restricted
product of characters $\nu_v=|\det|_v$ where $|\ \ |_v$ is the
$v$-adic norm and $\det$ is the reduced norm at the place $v$.
Let $m\in \n^*$ and $\rho\in DS_m$ be a cuspidal representation. If
$k\in\n^*$, then the induced representation to $G_{mk}(\aa)$ from
$\otimes_{i=0}^{k-1}(\nu^{\frac{k-1}{2}-i}\rho)$ has a unique
constituent (in the sense of \cite{LanCor}) $\pi$ which is a discrete series
(i.e. $\pi\in DS_{mk}$). We set then $\pi=MW(\rho,k)$. Discrete
series $\pi$ of groups $G_n(\aa)$, $n\in\n^*$, are all of this type,
$k$ and $\rho$ are determined by $\pi$. The discrete series $\pi$ is
cuspidal if $k=1$ and residual if $k>1$. These results are proved in
\cite{MW}. 

The proof of the following propositions and corollary is the same as in   \cite{Ba1}, 
once the local and global transfer are established without  condition on archimedean places. 
Firstly,  concerning cuspidal representations of $G'(\aa)$, we get:

\begin{prop}\label{cuspidal}
Let $m\in \n^*$ and let $\rho\in DS_m$ be a cuspidal representation.
Then

{\rm (a)} There exists $s_{\rho,{\mathcal D}}\in \n^*$ such that,
for $k\in \n^*$, $MW(\rho,k)$ is ${\mathcal D}$-compatible if and
only if $s_{\rho,{\mathcal D}}|k$. We have $s_{\rho,{\mathcal
D}}|d$.

{\rm (b)} ${\g}^{-1}(MW(\rho,s_{\rho,{\mathcal D}}))=\rho'\in DS'_{\frac{ms_{\rho,{\mathcal D}}}{d}}$
  is cuspidal. The map  $\g^{-1}$ sends cuspidal $\caD$-compatible representations to cuspidal representations.

{\rm (c)} Every cuspidal representation in
$DS'_{\frac{ms_{\rho,{\mathcal D}}}{d}}$ is obtained as in {\rm
(b)}.
\end{prop}

Let us call  {\it essentially cuspidal} representation the twist of
a cuspidal representation by a real power of $\nu$.
If $n_1,n_2,...,n_k$ are positive integers such that $\sum_{i=1}^k
n_i=n$, then the subgroup $L$ of $G'_n(\aa)$ of diagonal matrices by
blocks of sizes $n_1$, $n_2$,...,$n_k$ will be called {\it standard
Levi subgroup} of $G'_n(\aa)$. We identify then $L$ to
$\times_{i=1}^k G'_{n_i}(\aa)$. All the definitions extend in an obvious way  to $L$. 
The two statements in the following Proposition generalize respectively  \cite{MW} and  Theorem 4.4 in \cite{JS}.

\begin{prop}\label{classif}
{\rm (a)} Let $\rho'\in DS'_m$ be a cuspidal representation and let $k\in \n^*$. The induced representation from
$\otimes_{i=0}^{k-1}(\nu_{\rho'}^{\frac{k-1}{2}-i}\rho')$ has a
unique constituent $\pi'$ which is a discrete series, denoted by 
$\pi'=MW'(\rho',k)$. Every discrete series $\pi'$ of a group
$G'_n(\aa)$, $n\in\n^*$, is of this type, and $k$ and $\rho'$ are
determined by $\pi'$. The representation $\pi'$ is cuspidal if
$k=1$, and residual if $k>1$. If $\pi'=MW'(\rho',k)$, then
${\g}(\rho')=MW(\rho,s_{\rho,{\mathcal D}})$ if and only if
${\g}(\pi')=MW(\rho,ks_{\rho,{\mathcal D}})$.

{\rm (b)} Let $(L_i,\rho'_i)$, $i=1,2$, be such that $L_i$ is a
standard Levi subgroup of $G'_n(\aa)$ and $\rho'_i$ is an
essentially cuspidal representation of $L_i$ for $i=1,2$. Fix any
finite set of places $V'$ containing the infinite places and all the
finite places $v$ where $\rho'_{1,v}$ or $\rho'_{2,v}$ is ramified
(i.e. has no non-zero vector fixed under $K_v$). If, for all places
$v\notin V'$, the unramified subquotients of the induced
representations from $\rho'_{i,v}$ to $G'_n(\aa)$ are equal, then
 $(L_1,\rho'_1)$ and $(L_2,\rho'_2)$  are conjugate.
\end{prop}

We know by \cite{LanCor} that if $\pi'$ is an automorphic representation of
$G'_n$, then there exists  $(L,\rho')$ where $L$ is a
standard Levi subgroup of $G'_n$ and $\rho'$ is an essentially
cuspidal representation of $L$ such that $\pi'$ is a constituent of
the induced representation from $\rho'$ to $G'_n$. A corollary of
the point {\rm (b)} of the proposition is

\begin{cor}
 $(L,\rho')$ is unique up to conjugation.
\end{cor}

\subsection{Final comment and acknowledgment}

Let us say a word about the length of the paper which can be explained by our desire to give
complete proofs or/and references of  all the statements. For instance, the proof of $U(3)$ for 
$\GL(n,\bbH)$ in Section \ref{U3} is quite long in itself, and requires the material about Bruhat $G$-order
introduced in the previous section, not needed otherwise. We could have saved four or five pages by referring 
to \cite{Ta5} which gives the proof of $U(3)$ for $\GL(n,\bbR)$ and  $\GL(n,\bbC)$, but 
\cite{Ta5} is at the time still unpublished, and our arguments using Bruhat $G$-order could be used to simplify
the proofs in \cite{Ta5}. Our exposition is also intended for the reader who is interested in comparing the 
archimedean and non-archimedean theory, by making them as similar as possible. Our discussion of Vogan's 
classification in Section \ref{VOG} is also longer than  strictly needed, but we feel that it is important
that the relation between Vogan and Tadi\'c classifications is explained somewhere in some details.

We would like to thank D. Vogan for answering many questions concerning his work.

\section{Notation}

\subsection{Multisets}
Let  $X$ be a set.  We denote by $M(X)$ the set of functions from $X$ to $\bbN$ with finite
 support, and we consider  
 an element  $m \in M(X)$   as  a `set with multiplicities'. Such an element $m \in M(X)$  
will be typically 
 denoted by  
\[ m=(x_1,x_2,\ldots ,x_r)    \]
 It is a (non ordered) list of elements   $x_i$ in  $X$.

 The multiset  $M(X)$ is endowed with the structure of a monoid induced  from the one  on $\bbN$ : if  
$m=(x_1,,\ldots ,x_r)$, $n=(y_1,,\ldots ,y_s)$ are in  $M(X)$, we get 
\[m+n= (x_1,\ldots ,x_r,y_1,\ldots ,y_s) . \]

\bigskip

\subsection{Local fields and division algebras}

In the following, we will use the following notation :  $F$ is a local field,  
$|\, .\, |_F$ is the normalized absolute value on  $F$ and 
$A$ is a central division algebra over $F$  with  $\dim_F(A)=d^2$. 

If $F$ is archimedean, then either $F=\bbR$ and $A=\bbR$ or $A=\bbH$, the algebra of quaternions, or $F=A=\bbC$.

\bigskip

\subsection{$\GL$}
For $n \in \bbN^*$, we set $G_n=\GL(n,A)$ and  $G_0=\{1\}$.
We denote the reduced norm on $G_n$ by 
\[ RN\, : \,  G_n \rightarrow F^\times.  \] 
We set :
\[ \nu_n \,  :  G_n \rightarrow |RN(g)|_F.      \]
 When the value of $n$ is not relevant  to the discussion, we will simply put $G=G_n$ and $\nu=\nu_n$.

\begin{rmq} If  $A=F$, the reduced norm is just the  determinant.
\end{rmq}
The character $\nu$  of  $G$ is unramified and in fact the group of unramified characters of $G$ is  
\[ \caX(G)=\{ \nu^s, \, s \in \bbC  \}.\]

If $G$ is one of the groups $G_n$, or more generally, the group of rational points of any reductive 
algebraic connected group defined over $F$, we denote by  $\caM(G)$  the category of smooth representations of 
 $G$ (in the non archimedean case), or the category of    Harish-Chandra modules (in the archimedean case), 
with respect to a fixed maximal compact subgroup $K$ of $G$. For $\GL(n,\bbR)$, $\GL(n,\bbC)$ and $\GL(n,\bbH)$, 
these maximal compact subgroups are respectively chosen to be $\Or(n)$, $\U(n)$ and $\Sp(n)$, embedded
in the  standard way.  Then $\caR(G)$ denotes  the 
  Grothendieck group of the category of finite length representations in $\caM(G)$. This is the free $\bbZ$-module 
with basis  $\Irr(G)$, the  set of  equivalence classes of irreducible representations in  $\caM(G)$.
If  $\pi \in \caM(G)$, of finite length, we will again denote by $\pi$  its image in $\caR(G)$.  
When confusion may occur,
we will state precisely if we consider $\pi$ as a representation or as an element in $\caR(G)$.   

Set 
\[ \Irr_n=\Irr(G_n),\qquad  \Irr= \coprod_{n \in \bbN^*} \Irr_n,\quad \caR=\bigoplus_{n \in \bbN} \caR(G_n),. \]
If  $\tau \in \caM(G_n)$ or $\caR(G_n)$, we set  $\deg \tau=n$.

\subsection{Standard parabolic and Levi subgroups}

Let  $n \in \bbN$ and let  $\sum_{i=1}^s n_i=n$ be a partition of  $n$. The group  
\[ \prod_{i=1}^s G_{n_i} \]
is identified with the subgroup of  $G_n$ of bloc diagonal matrices of respective size 
 $n_1,\ldots, n_s$. Let us denote  $G_{(n_1,\ldots, n_s)}$ this subgroup and 
$P_{(n_1,\ldots, n_s)}$ (resp. $\bar P_{(n_1,\ldots, n_s)}$) the  parabolic subgroup of  $G_n$ containing 
  $G_{(n_1,\ldots, n_s)}$ and the  Borel subgroup of invertible upper triangular matrices 
(resp. lower triangular). The subgroup    $G_{(n_1,\ldots, n_s)}$ is a Levi factor of the  standard parabolic subgroup 
 $P_{(n_1,\ldots, n_s)}$.

In this setting, we denote by  $i_{(n_1,\ldots, n_s)}$  (resp.  $\underline i_{(n_1,\ldots, n_s)}$) 
the functor of normalized parabolic induction from $\caM(G_{(n_1,\ldots, n_s)})$ to  $\caM(G_n)$ with respect to 
 the parabolic subgroup  $P_{(n_1,\ldots, n_s)}$ (resp. $\bar P_{(n_1,\ldots, n_s)}$).
 
\bigskip

\begin{defi} Let  $\pi_1 \in \caM(G_{n_1})$ and  $\pi_2 \in \caM(G_{n_2})$, both of finite length. 
We can then form the induced representation : 
\[\pi_1\times \pi_2:= i_{(n_1,n_2)}(\pi_1 \otimes \pi_2).\]
We still  denote by  $\pi_1\times \pi_2$ the image  of $i_{n_1,n_2}(\pi_1 \otimes \pi_2)$ 
 in the Grothendieck group  $\caR_{n_1+n_2}$.  This extends linearly to a product 
\[\times: \caR \times \caR \rightarrow \caR. \]
\end{defi}

\bigskip 

\begin{rmq}
Again we warn the reader that it is important to know when we consider $\pi_1\times \pi_2$ as a representation 
or an element in $\caR$.   For instance 
$\pi_1 \times \pi_2=\pi_2 \times \pi_1$ in $\caR$ (see below), but  $i_{(n_1,n_2)}( \pi_1 \otimes \pi_2)$ is 
not isomorphic to 
$i_{(n_1,n_2)}( \pi_2 \otimes \pi_1)$ in general.
\end{rmq}

\begin{prop}
The ring   $(\caR,\times)$ is graded  commutative. Its identity is the unique element in $\Irr_0$. 
\end{prop}

\bigskip

\section{Langlands classification}

We recall how to  combine Langlands classification of $\Irr$ in terms of irreducible 
essentially tempered representations, and the fact that for the groups $G_n$, tempered representations
are induced fully from  irreducible square integrable modulo center representations to give a classification 
of $\mathbf{Irr}$ in terms of  irreducible essentially square integrable modulo center representations.

Let us denote respectively by  
\[D_n^u \subset \Irr_n, \qquad D_n \subset \Irr_n, \]
the set of equivalence classes of  irreducible,  square integrable modulo center
(respectively   essentially square integrable modulo center)  representations 
of $G_n$ and set  
\[ D^u= \coprod_{n \in \bbN^*} D_n^u, \qquad  D= \coprod_{n \in \bbN^*} D_n.  \]

Similarly, 
\[T^u_n \subset \Irr_n,\qquad T_n \subset \Irr_n, \]
denote respectively the sets of equivalence classes of  irreducible tempered  representations of $G_n$ and
 equivalence classes of  irreducible essentially  tempered  representations of $G_n$. Set   
\[ T^u= \coprod_{n \in \bbN^*} T^u_n, \qquad  T= \coprod_{n \in \bbN^*} T_n.  \]

For all  $\tau \in T$, there exists a unique  $e(\tau) \in \bbR$ and a unique  
$\tau^u \in T^u$     such that 
\[ \tau= \nu^{e(\tau)}  \tau^u. \]

\begin{thm} \label{Rtriv} Let  $d=(\delta_1,\ldots,\delta_l) \in M(D^u)$. Then  
\[ \delta_1\times \delta_2\times \ldots \times \delta_l \]
is irreducible, therefore in  $T^u$. 
This defines a one-to-one correspondence between  $M(D^u)$ and $T^u$.
\end{thm}

This is due to Jacquet and Zelevinsky in the case $A=F$ non archimedean (\cite{Jac} or \cite{Ze}). For a non archimedean
division algebra, this is established in \cite{DKV}. In the archimedean case, reducibility of induced from 
square integrable representations are well-understood in terms of $R$-groups (Knapp-Zuckermann \cite{KZ}), 
and for the groups $G_n$, the $R$-groups are trivial.

\begin{defi}
Let  $t =(\tau_1,\ldots,\tau_l)\in M(T)$. We say that  $t$ is written in a standard order if 
\[ e(\tau_1) \geq  \ldots \geq e(\tau_l). \]
\end{defi}
\bigskip

\begin{thm}\label{Lgclass}
Let  $d=(d_1,\ldots ,d_l) \in  M(D)$  written in a standard order, {\sl i.e.}
\[ e(d_1) \geq  e(d_2)> \ldots   \geq e(d_l). \]

Then :  

--- $(i)$ the representation :
\[ \lambda(d)=d_1 \times \ldots  \times d_l \]
has a  unique  irreducible quotient $\Lg(d)$, appearing with multiplicity one in  
a Jordan-H\"older sequence of  $\lambda(d)$. 
It is also the unique subrepresentation of 
\[ d_l\times d_{l-1} \times \ldots \times d_2 \times d_1. \]

--- $(ii)$ Up to a multiplicative scalar, there is an unique intertwining operator   
\[ J : d_1 \times \ldots  \times d_l \longrightarrow  d_l \times \ldots  \times d_1.  \]
We have then  $\Lg(d) \simeq \lambda(d)/\ker J\simeq \im J$.

--- $(iii)$  The map  
\[ d \mapsto \Lg(d) \]
is a bijection between  $M(D)$ and $\Irr$. 
\end{thm}

For a proof in the non archimedean case, the reader may consult \cite{Re}.

\bigskip 

Representations of the form  $\lambda(d)=d_1 \times \ldots  \times d_l$ with  
 $d=(d_1,\ldots ,d_l) \in  M(D)$  written in a standard order are called {\sl standard representations}.

\begin{rmq}
If $d$ is a multiset of representations in $\Irr$, we denote by $\deg d$ the sum of the degrees of representations 
in $d$. Let $M(D)_n$ be the subset of $M(D)$ of multisets  of degree $n$. Then the theorem gives a one-to-one
 correspondence between  $M(D)_n$ and $\Irr_n$.  
\end{rmq}

\bigskip

\begin{prop}
The ring  $R$  is isomorphic to  $\bbZ[D]$, the ring of polynomials in $X_d$,  
 $d \in D$  with coefficients in $\bbZ$, {\sl i.e.} $\{[\lambda(d)]\}_{d \in D}$   is a  $\bbZ$-basis of
 $\caR$.
\end{prop}

See \cite{Ze}, Prop. 8.5 for a proof.

\bigskip

We give some easy consequences  of the proposition : 
\begin{cor} $(i)$ The ring   $R$  is a factorial domain.

--- $(ii)$ If  $\delta \in D$, $[\delta]$ is prime  $\caR$.

--- $(iii)$ If  $\pi \in \caR$  is homogeneous and $\pi=\sigma_1 \times \sigma_2$ 
in  $ \caR$ ,  the  $\sigma_1$ and  $\sigma_2$ are homogeneous.

--- $(iv)$ The group of invertible elements in  $ \caR$  is $\{ \pm \Irr_0 \}.$
\end{cor}

\bigskip

\section{Jacquet-Langlands correspondence}\label{JL1}

In this section, we fix  a central division algebra $A$ of dimension $d^2$ over the local field $F$.
We recall the Jacquet-Langlands correspondence between $\GL(n,A)$ and $\GL(nd,F)$. Since we need simultaneously
both $F$ and $A$ in the notation, we set $G_n^A$, $G_n^F$ respectively for $\GL(n,A)$, $\GL(n,F)$, and similarly with 
other notation {\sl e.g.} $\caR(G^A_n)$ or  $\caR(G^F_n)$, $D_n^A$ or $D_n^F$, etc.

There is a standard way of defining the determinant and the
characteristic polynomial for elements of $G_n$, in spite of $A$
being non commutative (see for example \cite{Pi} Section  16), and the reduced norm $RN$ introduced above
is just given by  the constant term of the characteristic polynomial. If $g\in
G_n$, then the characteristic polynomial of $g$ has coefficients
in $F$, it is monic and has degree $nd$. If $g\in G_n$ for some $n$, we say 
$g$ is {\sl regular semisimple} if the characteristic polynomial of
$g$ has distinct roots in an algebraic closure of $F$.

If $\pi\in \caR(G_n)$, then we let $\Theta_\pi$ denote the function character of $\pi$, 
as a locally constant map, stable under conjugation, defined on the set
of regular semisimple elements of $G_n$.

 We say that $g' \in G_n^A$ {\it corresponds}
to $g \in G_{nd}^F$ if $g$ and $g'$ are
regular semisimple and have the same characteristic polynomial, and we   write then $g'\leftrightarrow  g$.
Notice that if  $g'\leftrightarrow  g$ and if $g_1'$ and $g_1$ are respectively conjugate to $g'$ and $g$, then 
$g'_1\leftrightarrow  g_1$. Said otherwise, it means that  $\leftrightarrow $ is really a correspondence between 
conjugacy classes.

\begin{thm}
There is a unique bijection $\mathbf{C}:D_{nd}^F \rightarrow  D_{n}^A$ such that for all 
$\pi\in D^F_{nd}$ we have
$$\Theta_\pi(g)=(-1)^{nd-n}\Theta_{\mathbf{C}(\pi)}(g')$$
for all $g\in G^F_{nd}$ and $g'\in G^A_n$ such that $g'\leftrightarrow  g$.
\end{thm}

For the proof, see \cite{DKV} if the characteristic of the base field $F$ is zero 
and \cite{Ba2} for the non zero characteristic case. In the archimedean case, see sections \ref{AR} and \ref{AH}, 
and Remark \ref{JLarch} for more details about this correspondence (\cite{JL},\cite{DKV}).

We identify the centers of $G^F_{nd}$ and $G^A_n$ via the canonical isomorphisms with $F^\times$.  Then the
correspondence $\mathbf{C}$ preserves central characters so in particular $\sigma$ is unitary if and only if 
$\mathbf{C}(\sigma) $ is.

The correspondence $\mathbf{C}$ may be extended in a natural way to a correspondence $\LJ$ 
between  Grothendieck groups : 

- If $\sigma \in D^F_{nd}$, viewed as an element in $\caR(G^F_{nd})$,  we set  $$\JL(\sigma)=(-1)^{nd-n} 
\ccc(\sigma),$$ 
viewed as an element in  $\caR(G^A_{n})$. 

- If $\sigma \in D^F_{r}$, where $r$ is not divisible by $d$, we set $\JL(\sigma)=0$. 

- Since $\caR^F$ is a polynomial algebra in the variables $d \in D^F$, one can extend $\JL$ in a unique way 
to an algebra morphism between $\caR^F$ and $\caR^A$. It is clear that $\JL$ is surjective.

The fact that $\JL$ is a ring morphism means that ``it commutes with parabolic induction''. 
Let us describe how to compute (theoretically) $\JL(\pi)$,  $\pi \in \caR^F$. 
Since  $\{ \lambda(a)\}_{a\in M(D^F)}$ is a basis of $\caR^F$, we first write $\pi$ in this basis as 
\[ \pi=\sum_{a\in M(D^F)} M(a,\pi) \lambda(a),   \]
with $M(a,\pi) \in \bbZ$  (see Section \ref{BrGorder}).
Since $\JL$ is linear, 
\[ \JL(\pi)=\sum_{a\in M(D^F)} M(a,\pi) \; \JL(\lambda(a)),    \]
so it remains to describe $\JL(\lambda(a))$. If $a=(d_1,\ldots d_k)$, then 
\[\lambda(a)=d_1\times \ldots \times d_k \]
(since we consider $\lambda(a)$ as an element in $\caR^F$, the order of the $d_j$ is not important).
Since $\JL$ is an algebra morphism 
\[\JL(\lambda(a))= \JL(d_1)\times \ldots \times \JL(d_k).  \]
If  $d$ does not divide one of the  $\deg d_i$, this is $0$, and   
if  $d$  divides  all the  $\deg d_i$, setting  $\sum_i \deg d_i=nd$, we get 
\[ \prod_{i=1}^k (-1)^{\deg d_i-\deg d_i/d}   
 \ccc(d_1)  \times \ldots \times \ccc(d_k)= (-1)^{nd-n}  \ccc(d_1)  \times \ldots \times \ccc(d_k).\]

\section{Support  and infinitesimal character}

The goal of this section is again to introduce some notation and to recall well known results, but we want 
to adopt a uniform terminology for archimedean and non archimedean case. In the non archimedean case, some authors, 
by analogy with the archimedean  case, call  `infinitesimal character' the cuspidal support of a representation 
(a multiset of irreducible supercuspidal representations).
We take the opposite view of considering infinitesimal characters in the archimedean case as multisets of 
complex numbers.

\subsection{Non archimedean case }\label{nonarchsupp}

We start with the case  $F$  non archimedean. We  denote  by $C$ (resp. $C^u$) the subset of  $\Irr$ 
of  supercuspidal representations  (resp. unitary supercuspidal, {\sl i.e.} such that  $e(\rho)=0$).

For all  $\pi \in \Irr$, there exist  $\rho_1,\ldots,\rho_n \in C$ such that  
$\pi$ is a subquotient of $\rho_1\times \rho_2\times \ldots \times \rho_n$. 
The multiset  $(\rho_1,\ldots,\rho_n) \in M(C)$ is uniquely determined by  $\pi$, and we denote it by  
$\supp(\pi)$. It is called the cuspidal support of $\pi$.  
When  $\pi$ is a finite length representation whose  irreducible subquotients  have same 
cuspidal support, we denote it by  $\supp (\pi)$.
If  $\tau=\pi_1\times \pi_2$, with  $\pi_1,\pi_2 \in \Irr$ we have  
\begin{equation}\label{S1} \supp (\tau)=\supp (\pi_1)+\supp(\pi_2)\end{equation}

For all  $\omega \in M(C)$, denote by  $\Irr_\omega$ the set of  $\pi \in \Irr$ whose cuspidal  support 
is  $\omega$. We obtain a decomposition : 
 \begin{equation} \label{S2} \Irr=\coprod_{\omega \in M(C)} \Irr_\omega. \end{equation}  
Set 
\[  \caR_\omega= \bigoplus_{\pi \in \Irr_\omega} \bbZ\,  \pi.\]
Then  
\begin{equation}\label{S3}  \caR=\bigoplus_{ \omega \in M(C)}  \caR_\omega \end{equation}
is a graduation of  $\caR$ by  $M(C)$. 

We recall the following well known result. 
\begin{prop}\label{finiNA} Let  $\omega \in M(C)$. Then   $\Irr_\omega$ is  finite.
\end{prop}

\bigskip 

\subsection{Archimedean case}

  Denote by  $\frg_n$ the  complexification of the  Lie algebra of  
$G_n$, $ \frU_n=\frU(\frg_n)$ its  enveloping algebra, and  $\frZ_n$ the center of the latter.
Let  $\frh_n$ be a Cartan subalgebra of  $\frg_n$, and  $W_n=W(\frg_n,\frh_n)$ its  Weyl group. 
Harish-Chandra has defined an algebra isomorphism from $\frZ_n$ to the Weyl group invariants in the 
symmetric algebra over $\frh_n$ :
\[ \mathrm{HC}_n \, : \,  \frZ_n \longrightarrow S(\frh_n)^{W_n}.  \]
Using this isomorphism, every character of $\frZ_n$ ({\sl i.e.} a morphism of algebra with unit
$ \frZ_n \rightarrow \bbC$) is identified with a character of  $S(\frh_n)^{W_n}$. Such characters are  
 given by orbits of  $W_n$ in  $\frh_n^*$, by evaluation at a point of the orbit.

A representation  (recall that in the archimedean case, this means a Harish-Chandra module) admits an 
 infinitesimal character if the center of the enveloping algebra acts on it  by scalars. Irreducible representations
admit infinitesimal character.  For all $\lambda \in \frh_n^*$, let us denote by  $\Irr_\lambda$ the set of  
$\pi \in \Irr$ whose  infinitesimal character is   given by   $\lambda$.

We are now going to identify  infinitesimal characters with multisets of complex numbers.

--- $A=\bbR$. In this case, $\frg_n=M_n(\bbC)$ and we can choose  $\frh_n$ to be the space of diagonal matrices,
identified with $\bbC^n$. Its dual space is also identified with  $\bbC^n$ by the canonical duality
\[\bbC^n \times \bbC^n\rightarrow \bbC, \quad ((x_1,\ldots ,x_n),(y_1,\ldots y_n))\mapsto \sum_{i=1}^n x_iy_i.  \]
The Weyl group $W_n$ is then identified with the  symmetric group  $\frS_n$, acting on  $\bbC^n$ by 
permuting coordinates.
 Thus, an  infinitesimal character for  $G_n$ is given by a multiset of   $n$  complex numbers.

--- $A=\bbC$.   In this case, $\frg_n=M_n(\bbC)\oplus M_n(\bbC) $, and we can choose   $\frh_n$
to be the space of couples of diagonal matrices, identified with  $\bbC^n \times \bbC^n$.    
 Its dual space is also identified with  $\bbC^n \times \bbC^n$ as above.
The Weyl group is then identified with $\frS_n \times \frS_n$, acting on  $\frh_n^*\simeq \bbC^n \times \bbC^n$
 by permuting coordinates. Thus, an infinitesimal character for  $G_n$  is given by a couple of multisets 
of  $n$   complex numbers.

--- $A=\bbH$. The group  $G_n$ is  a real form of  $\GL(2n,\bbC)$, so  $\frg_n=M_{2n}(\bbC)$.
The discussion is then the same as for  $F=\bbR$, with  $2n$ replacing  $n$. 

\bigskip 

 By analogy with the non archimedean case, we denote by $M(C)$ the set of  multisets 
(or couple of multisets if $A=\bbC$) described above.

\begin{defi} Let   $\omega \in M(C)$ be a  multiset 
(or a couple of  multisets of the same cardinality, if $A=\bbC$) of  complex numbers.
  If  $\pi \in \Irr_n$, we set 
\[ \supp(\pi)=\omega \]
where $\omega \in M(C)$ is the multiset (or couple of multisets  of the same cardinality 
if $F=\bbC$) defined by the infinitesimal 
character of  $\pi$. We say that $\omega$ is the support of $\pi$. 
 When $\pi$ is a finite length  representation whose  subquotients have all same support, we denote it by
  $\supp (\pi)$.
If $\pi \in \mathbf{Irr}$, $\pi=\Lg(a)$ for  $a\in M(D)$, we set 
\[ \supp(a):=\supp(\pi).  \]
We denote by $M(D)_\omega$ the set of $a \in M(D)$ with support $\omega$. 
\end{defi}

\bigskip 

\begin{prop} The results of \ref{nonarchsupp} are valid in the archimedean case.
\end{prop}

By this, we mean (\ref{S1}),  (\ref{S2}), (\ref{S3})  and Prop. \ref{finiNA} above.

\bigskip

\section{Bruhat $G$-order} \label{BrGorder}

We continue with the notation of the previous section. In the sequel, we will use  a partial order $\leq$ on 
$M(D)$, called Bruhat   $G$-order,  obtained from partial orders 
on each $M(D)_\omega$, $\omega \in M(C)$ whose main properties are described in the following :

\begin{prop}\label{multmat}
Let  $a\in M(D)$. Then the decomposition of $\lambda(a)$ in the basis  $\{\Lg(b)\}_{b \in M(D)}$ of  $\caR$
is of the form 
\[ \lambda(a)=  \sum_{b\leq a} m(b,a)\; \Lg(b),   \]
where the $m(a,b)$ are  non negative integers.
The  decomposition  of $\Lg(a)$ in the basis  $\{\lambda(b)\}_{b \in M(D)}$ of  $\caR$
is of the form 
\[ \Lg(a)= \sum_{b\leq a} M(b,a)\; \lambda(b), \]
where the  $M(b,a)$ are integers.
In particular, all factors   $\Lg(b)$ (resp. $\lambda(b)$) appearing in  the decomposition of  
$ \lambda(a)$ (resp. $\Lg(a)$) have same support.
Furthermore,  $m(a,a)=M(a,a)=1$. 
\end{prop}

In the non archimedean case,  Bruhat $G$-order is described  by  Zelevinsky \cite{Ze} ($A=F$) and Tadi\'c \cite{Ta2}
 in terms of linked segments.
 On arbitrary real reductive  groups,  Bruhat $G$-order is  defined by  Vogan on a different sets of parameters,
 in terms of integral roots (see \cite{VogIC4}, def. 12.12). In all cases, 
 Bruhat $G$-order is constructed  by defining first  {\sl elementary operations}, starting from an 
element $a \in M(D)$ and  obtaining another  element  $a'\in M(D)$. This is written 
\[a' \prec a.\]
  Bruhat $G$-order is then  generated by  $\prec$\footnote{In the case  $A=\bbR$, the  situation 
is a little more complicated.}. Another important property of  Bruhat $G$-order is the following. 
 One can define on all  $M(D)_\omega$ a length function :
\[  l: \, M(D)_\omega \rightarrow \bbN  \]
 such that  if $b \leq a$, then  $l(b)\leq l(a)$,   if $b \leq a$ and   $l(b)= l(a)$ then  $b=a$ and 
finally if   $b \leq a$, and   $l(b)= l(a)-1$ then  $b \prec a$.  In particular, if  $b \prec a$, there is no  
$c \in M(D)_\omega$ such that  $b\leq c <a$ but  $b=c$.

We have then 
\begin{prop} \label{precA}
Let  $a,b\in M(D)_\omega$  such that   $b \prec a$.  Then   $m(b,a)\neq 0$ and $M(b,a)\neq 0$.
\end{prop}
\pr 
The first  assertion follows from the recursion formulas for Kazhdan-Lusztig-Vogan polynomials
in the  archimedean case \cite{VogIC3}. We even have  $m(b,a)=1$ in this case.
It is established by  Zelevinsky  \cite{Ze} or  Tadi\'c \cite{Ta2} in the non archimedean case,
 and  the  second  assertion  follows using  Prop.  \ref{multmat}. \qed 
\bigskip

\section{Unitary dual}
\label{UD}

\subsection{Representations $u(\delta,n)$ and  $\pi(\delta,n;\alpha)$}

Let  $\delta \in D$. Then  $\delta \times  \delta$ is irreducible. Indeed, if  $\delta \in D^u$, this is 
\ref{Rtriv}, and the general case follows by tensoring with an unramified character. 
Consider  $\delta \times \nu^\alpha \delta$, with  $\alpha>0$. 
There exists a smallest  $\alpha_0>0$ such that  $\delta \times \nu^{\alpha_0} \delta$ is reducible.

\begin{defi} \label{nudelta} Let  $\delta \in D$.  Set  $\nu_\delta=\nu ^{\alpha_0}$, where  
 $\alpha_0>0$ is the smallest real number  $\alpha>0$ such that  $\delta \times \nu^\alpha \delta$ is reducible.
\end{defi}

For all  $\delta \in D$, and for all  $n \in \bbN^*$ we set 
\begin{equation}\label{adelta} a(\delta,n)=( \nu_\delta^{\frac{n-1}{2}}\delta,  \nu_\delta^{\frac{n-1}{2}-1}
\delta, \ldots,   
\nu_\delta^{-\frac{n-1}{2}} \delta    )\in M(D), \end{equation}
\begin{equation}\label{udelta} u(\delta,n)= \Lg( a(\delta,n)). \end{equation}  

For all  $\delta \in D$, for all  $n \in \bbN^*$, and for all  $\alpha \in \bbR$, set 
\begin{equation}\label{pidelta}\pi(\delta,n; \alpha)= 
\nu_\delta^\alpha  u(\delta,n) \times \nu_\delta^{-\alpha}  u(\delta,n). \end{equation}

\subsection{ Tadi\'c hypotheses  $U(0), \ldots, U(4)$ and classification of the unitary dual}

We recall  Tadi\'c's  classification of the unitary dual of the groups $G_n$. For a fixed division algebra $A$,
consider the following hypotheses : 

$U(0)$ :  if $\sigma, \tau \in \Irr^u$, then $\sigma \times \tau \in \Irr^u$.

$U(1)$ : if $\delta \in D^u$ and $n \in \bbN^*$, then   $u(\delta,n) \in \Irr^u$.  

$U(2)$  : if  $\delta \in D^u$, $n \in \bbN^*$ and  $\alpha \in ]0,1/2[$,  then   $\pi(\delta,n;\alpha)
 \in \Irr^u$.

$U(3)$ : if  $\delta \in D$,  $u(\delta,n)$ is prime in  $\caR$. 

$U(4)$ : if  $a,b \in M(D)$, then $L(a)\times L(b)$ contains  $L(a+b)$ as a subquotient.

\bigskip

Suppose Tadi\'c's hypotheses are satisfied for  $A$. We have then the following :

\begin{thm}\label{TaCla}  The set  $\Irr^u$ is endowed with the structure of a free commutative monoid, with product 
 $(\sigma,\tau)\mapsto \sigma \times \tau$ and with basis 
\[\caB=\{ u(\delta,n), \, 
\pi(\delta,n;\alpha)\, |\,  \delta \in D^u, n\in \bbN^*, \alpha \in ]0,1/2[ \;  \}.\]
More explicitly, if  $\pi_1,\ldots ,\pi_k \in \caB$, then  $\pi_1\times \ldots \times \pi_k \in \Irr^u$
and if  $\pi \in  \Irr^u$, there exists  $\pi_1,\ldots ,\pi_k \in \caB$, unique up to permutation, such that  
$\pi=\pi_1\times \ldots \times \pi_k $.
\end{thm}

This is proved in  \cite{Ta3}, prop 2.1. The proof is formal.

\bigskip 

Let us first notice  that  $U(4)$ is  a quite simple  consequence of  Langlands classification,
 established by Tadi\'c for all $A$ in \cite{Ta4} (the proof works also for archimedean $A$, see \cite{Ta5}).
  It is also easy to see that  
$U(2)$ can be deduced from  $U(0)$ and  $U(1)$  by the following simple principle : 
if  $(\pi_t)_{t \in I}$, is a family of hermitian representations in  $\caM(G)$,
where  $I$ is an open interval containing $0$, continuous in a sense that we won't make precise here, and if  
 $\pi_0$ is unitary and  irreducible, then  $\pi_t$ is unitary  on the largest interval  $J\subset I$
 containing  $0$ where  
$\pi_t$  is irreducible (the signature  of the hermitian form  can change only when crossing reducibility points).
 Representations $\pi(\delta,n;\alpha)$, $\alpha \in \bbR$ are  hermitian, 
\[\pi(\delta,n;0) =u(\delta,n)\times u(\delta,n) \]
 is  unitary and  irreducible  ($U(0)$ and  $U(1)$), and  $\pi(\delta,n;\alpha)$   is irreducible for 
$\alpha \in ]-\frac{1}{2},\frac{1}{2}[$. See \cite{Ta5} and the references given there  for details. 

For the remaining $U(0)$, $U(1)$ and $U(3)$, the situation is more complicated. 

--- $U(3)$ is proved by Tadi\'c in the non archimedean case in \cite{Ta1}, and for $A=\bbR,\bbC$ in \cite{Ta5}.
We give below the proof for $A=\bbH$, following Tadi\'c's ideas.

--- $U(1)$ is proved by   Tadi\'c in the non archimedean case in \cite{Ta1} for the field case $A=F$. 
The generalization to all division algebra over $F$ is given by the authors in \cite{BR1}, using unitarity of 
some distinguished representations closely related to the $u(\delta,n)$ established by the first 
named author in \cite{Ba3}
by  global methods. For $F=\bbC$,  $u(\delta,n)$ is a unitary character, so the statement is obvious. For $F=\bbR$,  
 $U(1)$ was first proved by Speh in  \cite{Sp} using  global method. It can also be proved using Vogan's results
on cohomological induction (see details below). Finally, for $A=\bbH$,  $U(1)$ can be established using again the 
general  results on cohomological induction, and the argument in \cite{BR1}. A more detailed discussion 
of the archimedean case is in section \ref{U1ARC}.
    
--- $U(0)$  is by far the most delicate point. For $A=F$ non archimedean, it is established by Bernstein in 
\cite{Be}, using reduction to the mirabolic subgroup. For $A=\bbR$ or $\bbC$, the same  approach can be used, 
but some serious technical difficulties remained unsolved until the paper of Baruch \cite{Bar}.
For $A$ a general non archimedean division algebra, $U(0)$ is established by V. S\'echerre \cite{secherre}
using his deep results
on  Bushnell-Kutzko's type theory for the groups $\GL(n,A)$, which give Hecke algebras isomorphisms and allow 
one to reduce the problem to the field case (the proof also  uses in a crucial way  Barbash-Moy results
 on unitarity for Hecke algebras representations \cite{BM}). In the case $A=\bbH$, there is to our knowledge no 
written references, but it is well-known to some experts that this can be deduced from Vogan's classification 
of the unitary dual of $G_n$ in the archimedean case (\cite{Vog2}). Vogan's classification is conceptually
 very different 
from Tadi\'c's classification. It has its own merits, but the final result is quite difficult to state  and 
to understand, since it uses sophisticated concepts and techniques of the theory of real reductive groups.
 So, for people interested mainly 
in applications, to automorphic forms for instance, Tadi\'c's classification is much more convenient. 
In the literature,  before Baruch's paper was published, one can often find the statement of Tadi\'c's 
classification, with reference to Vogan's
 paper \cite{Vog2} for the proof. It  might not be totally obvious for non experts   to derive
 Tadi\'c's classification  from Vogan's. We take this opportunity to explain in this paper
 (see \S \ref{VOG} below) some aspects of Vogan's classification's,
how it is related to  Tadi\'c's classification and  how to deduce  $U(0)$ from it. Of course, an independent proof of
$U(0)$ would be highly desirable in this case. It would be even better to have an uniform proof of $U(0)$ 
for all cases,  but for this, new ideas are clearly needed. 

--- all these results are true if the characteristic of $F$ is positive (as explained in \cite{bhls}).

\section{Classification of generic irreducible unitary representations}\label{GEN}

From the classification of the unitary dual of $\GL(n,\bbR)$ given above and the classification of 
irreducible generic  representations of a real reductive groups  (\cite{VoGK}, \cite{Ko}), we deduce 
the classification of  generic irreducible unitary representations  of $\GL(n,\bbR)$. Let us
 first recall that Vogan gives a classification of    'large' irreducible representations of a 
quasi-split real reductive group ({\sl i.e.} having maximal Gelfand-Kirillov
 dimension),  that Kostant shows that  such a group admits generic representations if and only if 
the group is quasi-split, and that ``generic'' is equivalent  to ``large''. Therefore, Vogan's result can 
be stated as follows : 
\begin{thm} 
Any generic  irreducible  representation of any quasisplit real reductive group is
irreducibly induced from a generic  limit of discrete series, and conversely, a representation which is 
irreducibly induced from a generic  limit of discrete series is generic.
\end{thm}

Let us notice that in the above theorem, one can replace ``limit of discrete series'' by ``essentially tempered'',
 because according to   \cite{KZ}, any tempered representation is fully induced from a limit of discrete series.
In the case of $\GL(n,\bbR)$ all  discrete series are generic, so by Theorem \ref{Rtriv}, all essentially tempered 
representations are generic.

Let us denote by $\mathbf{Irr}^u_{gen}$ the subset of  $\mathbf{Irr}^u$ consisting of generic representations.
We have then the following specialization of theorem \ref{TaCla}.

\begin{thm}\label{Gen}  The set  $\Irr^u_{gen}$ is endowed with the structure of a free commutative monoid, 
with product 
 $(\sigma,\tau)\mapsto \sigma \times \tau$ and with basis 
\[\caB_{gen}=\{ u(\delta,1), \, 
\pi(\delta,1;\alpha)\, |\,  \delta \in D^u, \alpha \in ]0,1/2[ \;  \}.\]
More explicitly, if  $\pi_1,\ldots ,\pi_k \in \caB_{gen}$, then  $\pi_1\times \ldots \times \pi_k \in \Irr^u_{gen}$
and if  $\pi \in  \Irr^u_{gen}$, there exists  $\pi_1,\ldots ,\pi_k \in \caB_{gen}$, unique up to permutation, such that  
$\pi=\pi_1\times \ldots \times \pi_k $.
\end{thm}

\section{Classification of discrete series : archimedean case}

In this section, we  describe explicitly square integrable modulo center irreducible representations 
of $G_n$ in the archimedean case. In the case $A=\bbH$, we give also details about 
 supports, Bruhat $G$-order... Since the Bruhat $G$-order is defined 
by Vogan on a set of parameters for irreducible representations consisting of (conjugacy classes of)
characters of Cartan subgroups, we also describe the bijections between the various sets of parameters.

\subsection{ $A=\bbC$}\label{AC}

There are  square integrable modulo center irreducible representations of $\GL(n,\bbC)$ only when $n=1$. Thus 
\[ D=D_1=\Irr_1. \]
An element  $\delta \in D$ is then a character 
\[ \delta: \GL(1,\bbC)\simeq \bbC^\times \rightarrow \bbC^\times   \]

 Let  $\delta \in D$. Then there exists a unique  $n \in \bbZ$ and a unique 
 $\beta \in \bbC$ such that  
\[ \delta(z)=   |z|^{2\beta}  \left(\frac{z}{|z|}\right)^n = |z|_\bbC^\beta \left(\frac{z}{|z|}\right)^n.    \]
Let  $x,y \in \bbC$ satisfying  
\[ \begin{cases}x+y&=2\beta\\
x-y&=n.
\end{cases}\]
 We set, with the above notation (and abusively writing a complex power of a complex number),  
\[ \delta(z)=\gamma(x,y) =z^x\bar z^y. \]

The following is well-known. 
\begin{prop} \label{nudC}
Let  $\delta=\gamma(x,y) \in D$ as above. Then  $\delta\times \nu^\alpha \delta$ 
is reducible for  $\alpha=1$ and irreducible for  
$0\leq \alpha<1$. Thus  $\nu_\delta=\nu$ ({\sl cf.} \ref{nudelta}). 
In the case of reducibility $\alpha=1$, we have in $\caR$: 
\[ \gamma(x,y) \times \gamma(x+1,y+1)= \Lg((\gamma(x,y), \gamma(x+1,y+1)))+ \gamma(x,y+1) \times \gamma(x+1,y) \]
\end{prop}

\bigskip

\subsection{ $A=\bbR$}\label{AR}

There are square integrable modulo center irreducible representations of  $\GL(n,\bbR)$ only when  $n=1,2$ :
\[ D=D_1 \coprod D_2=\Irr_1 \coprod D_2. \]

Let us start with the parametrization of $D_1$. 
An element  $\delta \in D_1$ is a character 
\[ \delta: \GL(1,\bbR)\simeq \bbR^\times \rightarrow \bbC^\times   \]

Let  $\delta \in D_1$. Then there exists a  unique $\epsilon \in \{0,1\}$ and a  
unique $\alpha \in \bbC$ such that 
\[ \delta(x)=   |x|^{\alpha}  \sgn(x)^\epsilon, \quad (x \in \bbR^\times).    \]
We set 
\[  \delta=\delta(\alpha,\epsilon).  \]

\bigskip

Let us now give a  parametrization of $D_2$. 
Let  $\delta_1,\delta_2 \in D_1$. Then  $\delta_1 \times \delta_2$ is  reducible if and only if there exists
 $p\in \bbZ^*$ such that  
\begin{equation*} \delta_1 \delta_2^{-1}(x)=  x^p   \sgn(x), \quad (x \in \bbR^\times)  \end{equation*}
If  $\delta_i=\delta(\alpha_i,\epsilon_i)$, we rewrite these conditions as 
\begin{equation}\label{redR}\alpha_1-\alpha_2=p, \quad \epsilon_1-\epsilon_2=p+1 \, \mod 2 \end{equation}

If  $\delta_1 \times \delta_2$ is  reducible, we have in  $R$,
\begin{equation}\label{reducR1} \delta_1 \times \delta_2= \Lg((\delta_1,\delta_2))+ \eta(\delta_1,\delta_2)  
    \end{equation}
where  $\eta(\delta_1,\delta_2) \in D_2$ and  $\Lg((\delta_1,\delta_2))$ is an irreducible 
 finite dimensional representation   (of dimension $|p|$ with the notation above). 

\begin{defi}
If  $\alpha_1,\alpha_2 \in \bbC$ satisfy  $\alpha_1-\alpha_2 \in \bbZ^*$, we set 
\begin{equation} \label{etaxy} \eta(\alpha_1,\alpha_2)=\eta(\delta_1,\delta_2) \end{equation}
where  $\delta_1(x)=  |x|^{\alpha_1}$ and $\delta_2(x)=|x|^{\alpha_2} \sgn(x)^{\alpha_1-\alpha_2+1}$. 
This define a surjective map from  
\[ \{ (\alpha_1,\alpha_2)\in \bbC^2\, | \, \alpha_1-\alpha_2 \in \bbZ^*  \} \]
to   $D_2$ and  
\[ \eta(\alpha_1,\alpha_2)=\eta(\alpha'_1,\alpha'_2)   \Leftrightarrow \{ \alpha_1,\alpha_2 \}= 
 \{ \alpha'_1,\alpha'_2 \}. \]

This gives a parametrization of $D_2$ by pairs of complex numbers  
$\alpha_1,\alpha_2$ satisfying  $\alpha_1-\alpha_2 \in \bbZ^*$.
\end{defi}

\begin{rmq}  The representation   $\eta(x,y) \in D_2$, $x,y \in \bbC$, $x-y \in \bbZ^*$  
is obtained from the character   $\gamma(x,y)$  of  $\bbC^\times$ by some appropriate functor of 
 cohomological induction. But, even when  $x=y$, the functor of cohomological induction  
maps  $\gamma(x,x)$  to an irreducible essentially tempered representation  of  $\GL(2,\bbR)$, namely the 
{\sl limit of discrete series} $\delta(x,0)\times \delta(x,1)$, which is  an irreducible principal series.
\end{rmq}

For that reason, we set for  $x \in \bbC$ : 
\begin{equation}\label{limsd}  \eta(x,x):= \delta(x,0)\times \delta(x,1)\in \Irr_2  \end{equation}

\bigskip

\begin{prop} \label{nudR}
Let  $\delta \in D$. Then  $\delta\times \nu^\alpha \delta$ is reducible for  $\alpha=1$ and irreducible for  
$0\leq \alpha<1$. Thus  $\nu_\delta=\nu$ ({\sl cf.} \ref{nudelta}). 
\end{prop}

This is also well-known.
Let us be more precise, by giving the composition series for $\delta\times \nu \delta$.   
We start with the case $\delta=\delta(\alpha,\epsilon)\in D_1$. Then we get from 
(\ref{reducR1}) that we have in $\caR$, 
\begin{equation}\label{reducD1}
\delta(\alpha,\epsilon)\times \delta(\alpha+1,\epsilon)= \Lg(\delta(\alpha,\epsilon), \delta(\alpha+1,\epsilon))
+ \eta(\alpha,\alpha+1) .
\end{equation}

In the case where  $\delta=\eta(x,y)\in D_2$, $x-y=r \in \bbN^*$, we get if   $r\neq 1$,  
\begin{align}\label{FoR} &\eta(x,y)\times \eta(x+1,y+1)=\\ 
\nonumber \Lg(\eta(x,y), &\eta(x+1,y+1)) + 
\eta((x,y+1)\times \eta(x+1,y)).\end{align}

If   $r=1$, the situation degenerates :
\begin{align*} &\eta(x,y)\times \eta(x+1,y+1)=\\ 
\nonumber \Lg(\eta(x,y), &\eta(x+1,y+1)) + 
\eta(x,y+1)\times \eta(x+1,y).\end{align*}

Recall that our convention is that   $$\eta(y+1,y+1)= \delta(y+1,0)\times \delta(y+1,1)$$ is a limit of discrete series, 
thus :
\begin{align}\label{FoRd} &\eta(y+1,y)\times \eta(y+2,y+1)=\\ 
\nonumber \Lg(\eta(y+1,y), &\eta(y+2,y+1)) + 
 \delta(y+1,0)\times \delta(y+1,1) \times \eta(y+2,y).\end{align}

\bigskip

\subsection{$A=\bbH$}\label{AH}

Let us identify quaternions and  $2\times 2$ matrices of the form  
\[ \left( \begin{array}{cc}  \alpha&\beta\\ -\bar \beta&\bar \alpha \end{array} \right), 
 \quad \alpha,\beta \in \bbC. \]
The reduced norm is given by 
\[ RN  \left( \begin{array}{cc}  \alpha&\beta\\ -\bar \beta&\bar \alpha \end{array} \right)
 =|\alpha|^2+|\beta|^2 .\]
The group of invertible elements $\bbH^\times$ contains  $\SU(2)$, the kernel of the reduced norm. 
Thus we have an exact sequence
\[1 \rightarrow \SU(2) \hookrightarrow \bbH^\times \stackrel{RN}{\longrightarrow}
\bbR^\times _+ \rightarrow 1,\]
and we can identify   $\bbH^\times$ with the direct product  $ \SU(2) \times \bbR^\times _+$.

The group  $\GL(n,\bbH)$ is a real form of  $\GL(2n, \bbC)$, its elements are matrices composed of  $2\times 2$  
quaternionic matrices described above. Complex conjugacy  on  $\GL(2n, \bbC)$ for this real form is given 
on the $2\times 2$   blocs  by 
 \[ \left( \begin{array}{cc}  \alpha&\beta\\ \gamma&\delta \end{array} \right) \mapsto 
 \left( \begin{array}{cc}  \bar \delta& -\bar \gamma \\ - \bar \beta & \bar \alpha \end{array} \right).\]

A maximal compact subgroup of  $ \GL(n,\bbH)$ is then  $$\mathbf{Sp}(n)\simeq \mathbf{U}(2n)\cap
\GL(n,\bbH).$$ Its rank is  $n$, the rank of  $\GL(n,\bbH)$ is  $2n$ and the split rank of the center is
 $1$. Thus there are square integrable modulo center representations only when  $n=1$. 

For  $n=1$,  $D_1=\Irr_1$, all irreducible representations of $\bbH^\times$ are  essentially  square 
integrable modulo center, since $\bbH^\times$ is compact modulo center. 
 Harish-Chandra's parametrization  in this case is as follows : irreducible representations of  $\bbH^\times$ 
are parametrized by some characters of a fundamental Cartan subgroup, here we choose 
\[ \bbC^\times \hookrightarrow \bbH^\times, \quad \alpha \mapsto \left( \begin{array}{cc}  \alpha&0\\ 0
&\bar \alpha \end{array} \right),  \]
which is connected.
Characters of  $\bbC^\times$ were described in  section $F=\bbC$. They are of the  form 
$\gamma(x,y)$, $x-y \in \bbZ$. An irreducible representation of  $\bbH^\times$ is then parametrized
by a couple of complex numbers $(x,y)$, such that  $x-y \in \bbZ$. The  couples 
 $(x,y)$ and   $(x',y')$ parametrize  the same representation if and only if the characters 
 $\gamma(x,y)$ and   $\gamma(x',y')$ are conjugate under the Weyl group, {\sl i.e.} 
if the multisets  $(x,y)$ and   $(x',y')$ are equal. Furthermore $\gamma(x,y)$
corresponds to an irreducible representation if and only if  
 $x \neq y$. Let us denote $\eta'(x,y)$ the representation 
parametrized by the  multiset $(x,y)$, $x-y \in \bbZ^*$. It is  obtained  from the character 
$\gamma(x,y)$ of the Cartan subgroup  $\bbC^\times$ by cohomological induction. 

\begin{rmq} As opposed to the case  $A=\bbR$, when we induced cohomologically the character  
 $\gamma(x,x)$ of the Cartan subgroup  $\bbC^\times$ to  $\bbH^\times$, we get  $0$ : 
there is no limits of discrete series.  Thus we set $\eta'(x,x)=0$. 
\end{rmq}

\bigskip

\begin{rmq}\label{JLarch}
 Jacquet-Langlands correspondence (see Section \ref{JL1}) between  representations of  $\GL(1,\bbH)=\bbH^\times$ and 
essentially square integrable modulo center irreducible  representations  $\GL(2,\bbR)$ is given by 
\[ \ccc(\eta(x,y))=\eta'(x,y), \quad x,y \in \bbC, x-y\in \bbZ^*. \]
The representations   $\eta(x,y)$ and  $\eta'(x,y)$ are obtained by cohomological induction from the same 
 character  $\gamma(x,y)$ of the Cartan subgroup  $\bbC^\times$ of  $\GL(2,\bbR)$ and  $\bbH^\times$. 
In the case $x=y$ the construction still respect the Jacquet-Langlands character relation since both sides are equal to zero.
\end{rmq}

More generally let us give now the parametrization of irreducible representations of  $\GL(n,\bbH)$ 
by conjugacy classes of characters of  Cartan subgroups. The group    $\GL(n,\bbH)$  
has only one conjugacy class of Cartan subgroups, a representative being  $T_n$, which consist of  
$2\times 2$ bloc diagonal matrices of the form  
 $\left( \begin{array}{cc}  \alpha&0\\ 0 & \bar \alpha  \end{array} \right) $.
Thus $T_n \simeq (\bbC^\times)^n$ is connected, $\frt_n= \mathrm{Lie}(T)\simeq \bbC^n$ and $(\frt_n)_\bbC\simeq
(\bbC\oplus \bbC)^n$.

 Let $\Lambda$ be a character of  $T_n$. Its differential 
\[ \lambda=d\Lambda \, : \,  \frt_n \rightarrow \mathbf{Lie}(\bbC^\times) \simeq \bbC, \]
 is a  $\bbR$-linear map, with complexification the  $\bbC$-linear map 
\[ \lambda=d\Lambda \, : \,  \frt_\bbC\simeq (\bbC\oplus \bbC)^n \rightarrow \mathbf{Lie}(\bbC^\times) \simeq \bbC. \]
Such a linear form is given by a $n$-tuple  of  couples $(\lambda_i,\mu_i)$ such that  
$\lambda_i-\mu_i \in \bbZ$. 

Since  $T_n$ is connected, a character  $\Lambda$ of  $T_n$ is determined by its differential. We write 
\[ \Lambda=\Lambda(\lambda_1,\mu_1,\ldots ,\lambda_n,\mu_n)= \Lambda((\lambda_i,\mu_i)_{1\leq i\leq n})\]
if its differential is given by the  $n$-tuple of  couples $(\lambda_i,\mu_i)$ such that  
$\lambda_i-\mu_i \in \bbZ$.

Let  $\caP$ the set of characters  $\Lambda=\Lambda((\lambda_i,\mu_i)_{1\leq i\leq n})$  of the Cartan subgroup $T_n$ 
such that  $\lambda_i-\mu_i \in \bbZ^*$.

Irreducible representations  of   $\GL(n,\bbH)$ are parametrized by  $\caP$, 
 two characters $\Lambda_1$ and  $\Lambda_2$ giving the same irreducible representations if 
and only if they are conjugate under  $W(\GL(2n,\bbC),T_n)$. 
This group is isomorphic to $\{\pm 1\}^n\times \frS_{n}$. Its action on 
$\frt_\bbC\simeq (\bbC\oplus \bbC)^n$ is as follows :
each factor $\{\pm 1\}$ acts inside the corresponding factor  $\bbC\oplus \bbC$ by permutation, 
and  $\frS_n$ acts by permuting the  $n$ factors $ \bbC\oplus \bbC$. Thus we see that irreducible   representations 
of  $\GL(n,\bbH)$ are parametrized by multisets of cardinality $n$ of  pairs of  complex numbers
$(\lambda_i,\mu_i)$ such that   $\lambda_i-\mu_i \in \bbZ^*$. Since such a pair $(\lambda_i,\mu_i)$ 
parametrizes the representation  $\eta'(\lambda_i,\mu_i)$, we recover the Langlands parametrization of  $\Irr$ by 
 $M(D)$.  Let us denote by  $\sim$ the equivalence relation on  $\caP$ given by the Weyl group action 
$W(\GL(2n,\bbC),T)$.  We have described one-to-one correspondences 
\[ \caP/\sim \;  \simeq \Irr_n \simeq M(D)_n \]

Recall that a support for  $\GL(n,\bbH)$ is a   
multiset of  $2n$  complex numbers, {\sl i.e.} an element of  the quotient of $\frt_\bbC^* \simeq (\bbC\oplus\bbC)^n\simeq
\bbC^{2n}$,  by the action of the  Weyl group   $W_\bbC \simeq \frS_{2n}$.  

\begin{defi}
The support of a character  $\Lambda=\Lambda((\lambda_i,\mu_i)_{1\leq i\leq n}) \in \caP$ is the multiset
\[(\lambda_1,\mu_1,\ldots ,\lambda_n,\mu_n).  \] 
It does not depend on the equivalence class of $\Lambda$ for  $\sim$. If  $\Lambda \in \caP$ 
parametrizes the irreducible representation  $\pi$, we have  $\supp(\Lambda)=\supp(\pi)$. 
\end{defi}

This describes explicitly the map 
\[ \caP \rightarrow M(C), \quad \Lambda \mapsto \supp(\Lambda)\]        
and its fibers : two parameters  
\[\Lambda_1((\lambda_i^1,\mu_i^1))\quad \text{ and  } \quad  
\Lambda_2((\lambda_i^2,\mu_i^2)),\]
 have same support if and only if the multisets
 \[ (\lambda_1^1, \ldots ,\lambda_n^1,\mu_1^1, \ldots ,\mu_n^1,) \quad \text{ and }\quad 
 (\lambda_1^2, \ldots ,\lambda_n^2, \mu_1^2, \ldots ,\mu_n^2)\] are equal.
We denote by  $\caP(\omega)$ the fiber at $\omega$.

Let us give now the description of the Bruhat $G$-order, in terms of integral roots.
We have the decomposition of Lie algebra : 
\[ \mathrm{Lie}(\GL(2n,\bbC))=(\frg_{2n})_\bbC= (\frt_n)_\bbC \bigoplus_{\alpha \in R} \frg_\bbC^\alpha      \]
where $R=\{ \pm(e_i-e_j), \, 1\leq i<j\leq 2n \}$ is the usual root system of type  $A_{2n-1}$.

The roots $\pm(e_{2j-1}-e_{2j})$, $j=1,\ldots, n$ are imaginary compact, thus 
\[ \sigma\cdot e_{2j-1} = e_{2j}, \quad   j=1,\ldots, n.\]
Other roots are complex: for all  $j,l$, $1\leq j\neq l \leq l$.
\[ \sigma\cdot   (e_{2j-1}-e_{2l-1})= e_{2j}-e_{2l}, \quad  \sigma\cdot   (e_{2j-1}-e_{2l})= e_{2j}-e_{2l-1}.
 \]

Let us fix a support  $\omega$ and let  
 $\Lambda$ be a character of  $T_n$ such that  $\supp(\Lambda)=\omega$, say $\Lambda=\Lambda
((\lambda_i, \mu_i)_{i={1,\ldots ,n}})$, $\lambda_i-\mu_i \in \bbZ^*$.
Notice  that $W_\bbC \simeq \frS_{2n}$ doesn't act on 
 $ \caP(\omega)$, since the condition 
\[ \lambda_i-\mu_i \in \bbZ\]
might not hold anymore after some permutation of the  $\lambda_i$. 

Denote by  $W_\Lambda$ the subgroup of  $W_\bbC$ consisting of elements  $w$ such that  
\[ w\cdot (\lambda_i,\mu_i)_i-(\lambda_i,\mu_i)_i \in (\bbZ \times \bbZ)^n .  \]
 Then  $W_\Lambda$ is the  Weyl group of the root system  $R_\Lambda$  of  {\sl integral roots} for $\Lambda$,
 a root   $\alpha=e_k-e_l$  in  $R$ being integral for  $\Lambda$ if, when writing 
\[ \lambda_1,\mu_1,\lambda_2,\mu_2,\ldots ,\lambda_n,\mu_n=\nu_1,\ldots ,\nu_{2n} \] 
then  $\nu_k-\nu_l \in \bbZ$. 

Suppose that the support $\omega$ is regular, {\sl i.e.} all the  
$\nu_i$, $1\leq i\leq 2n$ are distinct.
We choose as a positive root system  $R_\Lambda^+\subset R_\Lambda$ 
the roots  $e_k-e_l$  such that  $\nu_k-\nu_l >0$. This defines simple roots.

Let us state first a necessary and sufficient condition for reducibility of standard modules (for
regular support).

\begin{prop}
Let   $a=(\eta'(\lambda_i,\mu_i)_{i=1,\ldots,n}) \in M(D)_\omega$, parametrized by the character  
$\Lambda=\Lambda((\lambda_i, \mu_i)_{i={1,\ldots ,n}})$ of  $T_n$. Suppose that the  support 
\[\omega=(\lambda_1,\mu_1,\ldots, \lambda_n,\mu_n)\] is regular.
Then  $\lambda(a)$ is reducible if and only if there exists a simple root $e_k-e_l$ in  $R_\Lambda^+$, complex, 
such that, if $e_k-e_l=e_{2i-1}-e_{2j-1}$, or  $e_k-e_l=e_{2i}-e_{2j}$,  $i\neq j$, then  
\[ \lambda_i-\lambda_j>0 \quad \text{ and  } \quad  \mu_i-\mu_j>0, \]
and if $e_k-e_l=e_{2i-1}-e_{2j}$,  or  $e_k-e_l=e_{2i}-e_{2j-1}$,  $i\neq j$, then 
 \[ \lambda_i-\mu_j>0 \quad \text{ and  } \quad  \mu_i-\lambda_j>0. \]

When  $\omega$ is not  regular, we still have a necessary condition for reducibility: 
if  $\lambda(a)$is  reducible, then there exists a   root $e_k-e_l$ in  $R_\Lambda^+$, not necessarily simple, 
but still satisfying the condition above.
\end{prop}

See \cite{VogIC4}.

\begin{defi}\label{prec1} We still assume  $\omega \in M(C)$ to be regular, and suppose that 
  $\Lambda\in \caP(\omega)$ satisfies the reducibility criterion above for the simple integral complex root   
$e_k-e_l$. Write 
\[\Lambda=\Lambda((\lambda_1,\mu_1),\ldots, (\lambda_n,\mu_n))=\Lambda((\nu_1,\nu_2),\ldots, (\nu_{2n-1},\nu_{2n}))  \]
Let $\Lambda' \in \caP(\omega)$, obtained from  $\Lambda$ by exchanging   $\nu_k$ and $\nu_l$, and let 
$a' \in M(D)_\omega$  corresponding to $\Lambda'$. We say that  
$a'$ is obtained from  $a$ by an  {\sl elementary operation}, 
 and we write  $a'\prec a$.  The Bruhat  $G$-order  on  $M(D)_\omega$ is the partial order generated by  $\prec$.    
\end{defi}
 
\bigskip

Let us now deduce from the reducibility criterion above the invariant  $\nu_\delta$ attached  ({\sl cf.}
 definition \ref{nudelta}) to an essentially square integrable modulo center irreducible representation  
 $\delta=\eta'(x,y)$, $x,y \in \bbC$, $x-y \in \bbZ^*$.
We may suppose that $x-y =r> 0$, since  $\eta'(x,y)=\eta'(y,x)$.

\begin{prop}
With the previous notation, $\nu_\delta=\nu$ if  $r>1$ and   $\nu_\delta=\nu^2$  if $r=1$. 
Since  $r$ is the dimension of $\delta$, we see that $\nu_\delta=\nu$ except when $\delta$ is a one-dimensional 
representation of $\GL(1,\bbH)$.

\end{prop}

\pr We want to study the reducibility of 
\[ \pi= \eta'(y+r,y)\times \eta'(y+r+\alpha,y+\alpha)\]
 for $\alpha>0$. The support of this representation is regular if and only if 
 $y+r,y,y+r+\alpha,y+\alpha$ are distinct, but since  
\[ y+r+\alpha > y+\alpha>y, \quad    y+r+\alpha > y+r>y,    \]
the support is regular except when  $r=\alpha$.
The representation  $\pi$ is the standard representation attached to the character 
 $$\Lambda=\Lambda((y+r+\alpha,y+\alpha),(y+r,y)).$$

If   $\alpha \notin \bbZ$, the support is regular, 
all integral roots are imaginary compact for  $\Lambda$, and then  $\pi$
is irreducible.

If  $\alpha=1$ and  $r\neq 1$,   the support is regular, 
all the roots are integral for  $$\Lambda((y+r+1,y+1),(y+r,y)),$$
 and   $e_1-e_3$  is a complex root, simple in  
\[R_\Lambda^+=\{ e_1-e_3, e_1-e_2,e_1-e_4, e_3-e_2, e_3-e_4,e_2-e_4  \},\] 
satisfying the reducibility criterion, since  
\[ (\sigma\cdot (e_1-e_3) )(y+r+1,y+1,y+r,y) = (e_2-e_4) (y+r+1,y+1,y+r,y) = 1>0. \]
The only smaller element that  $\Lambda$ in the Bruhat $G$-order is  
$$\Lambda'=\Lambda((y+r,y+1),(y+r+1,y)),$$ 
and we get 
\begin{align}\label{Fo} &\eta'(y+r,y)\times \eta'(y+r+1,y+1)=\\ 
\nonumber \Lg(\eta'(y+r,y), &\eta'(y+r+1,y+1)) + 
\eta'(y+r,y+1)\times \eta'(y+r+1,y).\end{align}

If  $\alpha=1$ and  $r=1$, the support is singular. Applying Zuckerman translation functors (see \cite{KV} 
for instance),  we get  
\begin{align*} &\eta'(y+2,y)\times \eta'(y+3,y+1)=\\ 
\nonumber \Lg(\eta'(y+2,y), &\eta'(y+3,y+1)) + 
\eta'(y+1,y+1)\times \eta'(y+2,y).\end{align*}

But, according to our convention, $\eta'(y+1,y+1)=0$ (this is really what we get applying
 translation functor to the wall), thus 
\begin{equation*} \eta'(y+1,y)\times \eta'(y+2,y+1)= \Lg(\eta'(y+1,y), \eta'(y+2,y+1)).\end{equation*}
is irreducible. 

The next possibility of reducibility for $r=1$ is then   $\alpha=2$, 
but then the support is regular and we see as above that 
$\pi$ is reducible, more precisely 
\begin{align} \label{redk1} &\eta'(y+3,y+2)\times \eta'(y+1,y)=\\
\nonumber \Lg(\eta'(y+3,y+2), &\eta'(y+1,y)) + 
\eta'(y+2,y+1)\times \eta'(y+3,y).\end{align}
\qed

\section{$U(3)$  for $A=\bbH$}\label{U3}

We follow  Tadi\'c \cite{Ta5} who gives a proof of $U(3)$ for $A=\bbC, \bbR$ to deal with the case $A=\bbH$.

\begin{thm} Let  $\delta=\eta'(y+r,y) \in D$ and  $y \in \bbC$ and $r \in \bbN^*$ and let $n \in \bbN^*$. 
Then  $u(\delta,n)$ is a prime in the ring   $\caR$. 
\end{thm}

\pr  We know that  $\delta$ is prime in  $\caR$, thus we start with  $n \geq 2$. 
Let us first deal with the case  $r=1$. Then $\nu_\delta=2$ and 
\[ a_0= a(\delta,n)=(\nu_\delta^{\frac{n-1}{2}}\delta,\nu_\delta^{\frac{n-1}{2}-1}\delta,\ldots ,
\nu_\delta^{-\frac{n-1}{2}}\delta ) =\]
\begin{small}
\[(\eta'(y+n, y+n-1), \eta'(y+n-2,y+n-3), \ldots 
 \eta'(y-n+2,y-n+1)) . \]\end{small}

Set  $a_0=a(\delta,n)=(X_1,\ldots,X_n)$ with  
\[ X_i= \gamma\left(y+n+2-2i, y+n+1-2i\right), \quad i=1,\ldots,n.   \]

\begin{rmq}\label{regu} The support of  $u(\delta,n)$ is the multiset 
\[ (y+n+2-2i, y+n+1-2i)_{i=1,\ldots n}.  \]
This support is regular.
\end{rmq}

Suppose that  $u(\delta,n)$ is not prime in  $\caR$. 
Then there exists polynomials $P$ and $Q$ in the variables $d \in D$, non invertible, such that 
$u(\delta,n)=PQ$. Since  $u(\delta,n)$ is homogeneous in  $\caR$ for the natural graduation, the same holds for
$P$ and $Q$.

Let us write  
\begin{equation}\label{mPQ}
P=\sum_{c \in M(D)} m(c,P)\lambda(c), \quad Q=\sum_{d \in M(D)} m(d,Q)\lambda(d).
\end{equation}

Set  $S_P=\{ a \in M(D)  |  m(a,P)\neq 0 \}$, $S_Q=\{ a \in M(D)  |  m(a,Q)\neq 0 \}$.
We get  
\[\Lg(a_0)=X_1\times X_2\ldots\times  X_n +\sum_{a \in M(D), \, a<a_0}  M(a,a_0) \; \lambda(a).   \]
Thus there exists $c_0 \in S_P$ and $d_0 \in S_Q$ such that  
\[ c_0+d_0=a_0=(X_1,\ldots ,X_n). \] 
  Since  $\deg P>0$, $\deg Q>0$,  $c_0$ and  $d_0$ are not empty. Furthermore, we see that the polynomials  
$P$ and $Q$ are not  constant
Denote by  $S_1$ the set of $X_i$ such that  $X_i \in c_0$ and by  $S_2$ the set of  $X_i$ such that  $X_i \in d_0$.
We get a partition of the $X_i$'s in two non empty disjoint sets. Thus we can find 
$1\leq i\leq n-1$ such that 
\[\{ X_i,X_{i+1}\} \not\subset S_j, \quad j=1,2, \] and without any loss of generality, we may suppose that 
$X_i \in S_1$, $X_{i+1} \in S_2$. Furthermore, we have
\[ |S_1|=\deg P, \quad |S_2|=\deg Q, \quad \deg P+\deg Q=n \]

We saw above  (\ref{redk1})  that  $X_i\times X_{i+1}$ is reducible, more precisely 
\[ X_i\times X_{i+1} =\Lg(X_i,X_{i+1})+Y_i\times Y_{i+1}  \]
where  $Y_i=\eta'(y+n+2-2i,y+n-1-2i)$, $Y_{i+1}=\eta'(y+n+1-2i,y+n-2i)$.

We have  $  a_1:=(Y_i,Y_{i+1})\prec (X_i, X_{i+1})$. 
 Set  \[ a_{i,i+1}= a_1+   (X_1,\ldots,X_{i-1},X_{i+2},\ldots   ,X_n). \]
Then  $a_{i,i+1} \prec a_0$ and thus $M(a_1,a_0)\neq 0$ by Prop. \ref{precA}.  
Therefore,  there exists  $c_1 \in S_P, \, d_1 \in S_Q$, non empty and  such that  
\[c_1+d_1=a_{i,i+1}.\]
We suppose now that  $Y_i$ divide $\lambda(c_1)$  in $\caR$, the case where  $Y_i$ divide $\lambda(d_1)$ 
being similar.

Suppose that also $\lambda(Y_{i+1})$ divide  $\lambda(c_1)$. We get a partition of the 
$X_j$, $j\neq i,i+1$ in two non empty sets $S'_1$ and  $S'_2$, such that  
\[c_1=\{ X_j, j\in S'_1 \}+Y_i+ Y_{i+1}, \quad    d_1=\{ X_j, j\in S'_2 \}.      \]  
The polynomials  $P$ and  $Q$ being homogeneous, we get  
\[\deg (P)=|S'_1|+2, \quad \deg(Q)=|S'_2|.\]
We see that  $X_{i+1}\notin T:=S_1 \cup S'_2$, thus  $ \{X_1,\ldots ,X_n\} \not\subset T $. 
For  $r \in \caR$, let us denote by  $\deg_T (r)$ the degree of  $r$ in the variables  $X_j \in T$. 
We get  $\deg_T P\geq |S_1|=\deg P$, $\deg_T Q\geq |S'_2|=\deg Q$, thus $\deg_T (\Lg(a_0))\geq n$. But the total degree 
of $\Lg(a_0)$ being  $n$, we get    $\deg_T  (\Lg(a_0))=n$. The expression of   $\Lg(a_0)$ in the basis  
$\lambda(b)$, $b\leq a_0$  shows that we can find  $b_0 \in M(D)$, 
$M(b_0,a_0)\neq 0$,  $\deg(b_0)=n$ and   $\deg_T  \lambda(b_0)=n$.   
 Furthermore $\lambda(b_0)$ can be written
\[  \lambda(b_0)= X_1^{\alpha_1}X_2^{\alpha_2}\ldots X_n^{\alpha_n}, \quad \alpha_j \in \bbN, \alpha_1+\cdots+\alpha_n=n. 
   \]
Since  $T\neq \{X_1,\ldots ,X_n\}$, there exists  $j$ such that  $\alpha_j > 1$. 
But then  $X_j$ appears with multiplicity at least two in  $b_0$. Since  
$\supp (b_0)=\supp(a_0)$ is regular, we get a  contradiction.

Suppose now that  $\lambda(Y_{i+1})$ doesn't divide    $\lambda(c_1)$. 
We get a partition of the 
$X_j$, $j\neq i,i+1$ in two non empty sets  $S'_1$ and $S'_2$, such that  
\[c_1=\{ X_j, j\in S'_1 \}+Y_i, \quad    d_1=\{ X_j, j\in S'_2 \}+ Y_{i+1}.      \]  
We set now  $T=S'_1 \cup S_2$, and we see that  $X_{i+1}$ doesn't belong to $T$ , thus  
$ \{X_1,\ldots ,X_n\} \not \subset T$.  For  $r \in \caR$,  denote by $\deg_T (r)$ the degree of $r$ in 
the variables   $X_j \in T$ and $Y_i$. As above, we get  : $\deg_T(\Lg(a_0)) =n$, there exists 
$b_0 \in M(D)$,  $M(b_0,a_0)\neq 0$, $\deg(b_0)=n$ and    $\deg_T  (\lambda(b_0))=n$. 
We can write 
\[  \lambda(b_0)= X_1^{\alpha_1}X_2^{\alpha_2}\ldots X_n^{\alpha_n}Y_i^\alpha, \quad
 \alpha_j \in \bbN, \alpha_1+\cdots+\alpha_n+\alpha=n.    \]
Since  $\{X_1,\ldots ,X_n\} \not \subset T$,  there exists  $j$ such that  $\alpha_j=0$. If  $\alpha=0$,
 we get a contradiction as above.
Thus  $\alpha\geq 1$, but since multiplicities in  $\supp(a_0)$ are at most  $1$, we get  $\alpha=1$,
$\alpha_j=1$ if $j\neq i+1,i$, $\alpha_i=\alpha_{i+1}=0$ and we still get a contradiction. This finishes the case 
$r=1$.

Let us deal with briefly the case  $r>1$. Then $\nu_\delta=\nu$ and 
\[  a(\delta,n)=(\nu^{\frac{n-1}{2}}\delta,\nu^{\frac{n-1}{2}-1}\delta,\ldots ,\nu^{-\frac{n-1}{2}}\delta ) =\]
\begin{Small}
\[(\eta'(x+\frac{n-1}{2},y+\frac{n-1}{2}),\eta'(x+\frac{n-1}{2}-1,y+\frac{n-1}{2}-1), \ldots 
 \eta'(x-\frac{n-1}{2},y-\frac{n-1}{2}) ). \]\end{Small}

Set  $a_0=a(\delta,n)=(X_1,\ldots,X_n)$ with  
\[ X_i= \gamma\left(y+r+\frac{n-1}{2}+1-i, y+ \frac{n-1}{2}+1-i    \right), \quad i=1,\ldots,n.   \]
We  proceed as above, using now formula  (\ref{Fo})
for the reducibility of  $\lambda(X_i,X_{i+1})$ : 
\[  \lambda(X_i,X_{i+1})=\Lg(X_i,X_{i+1})+\lambda(Y_i,Y_{i+1})  \]
where  $Y_i=\eta'(y+r+\frac{n-1}{2}+1-i , y+ \frac{n-1}{2}-i)$ and     
$Y_{i+1}=\eta'(y+r+\frac{n-1}{2}-i , y+ \frac{n-1}{2}+1-i)$.
In all cases, we get  contradictions by inspecting   multiplicities in the 
 support. We leave the details to the reader. \qed

\bigskip

\section{$U(1)$ :  archimedean case} \label{U1ARC}

We recall briefly the arguments for $A=\bbC$ and $\bbR$, even if it is well known and done elsewhere, because 
we will need the notation anyway.  We give the complete argument when $A=\bbH$. 

\subsection{$A=\bbC$} \label{U1C}

This case is easy because  for $\gamma=\gamma(x,y)$, $x,y \in \bbC$, $x-y \in \bbZ$, a character of $\bbC^\times$,  
we have 
\[ u(\gamma,n) = \gamma \circ  \det .  \]
Representations   $u(\gamma,n)$ are thus  $1$-dimensional representations  of  $\GL(n,\bbC)$. Furthermore, 
if  $\gamma$ is unitary  ({\sl i.e.} $\Re e (x+y)=0$) then  $u(\gamma,n)$ is unitary.

\subsection{$A=\bbR$}\label{U1R}

There are two cases to consider. The first is  $\delta \in D_1$, $\delta=\delta(\alpha,\epsilon)$, 
$\alpha \in \bbC$, $\epsilon \in \{0,1\}$. This case is similar to the case  $A=\bbC$ above, since  
\[ u(\delta,n) = \delta \circ  \det .  \]
Representations  $u(\delta,n)$ are  $1$-dimensional representations  of   $\GL(n,\bbR)$. Furthermore, 
if  $\delta$ is unitary  ({\sl ie.} $\Re e (\alpha)=0$) then  $u(\delta,n)$ is unitary. 

The second case is  $\delta=\eta(x,y)\in D_2$, $x,y \in \bbC$, $x-y=r \in \bbN^*$. 
We have already mentioned without giving any details that  $\eta(x,y)$ 
is obtained by  cohomological induction from the character   
 $\gamma(x,y)$ of the Cartan subgroup  $\bbC^\times$ of $\GL(2,\bbR)$.
Let us be now more precise. Cohomological induction functors considered here are normalized as in \cite{KV}, (11.150b): 
if  $(\frg_\bbC,K)$ is a reductive pair associated to a real reductive group  $G$,  if 
 $\frqqq_\bbC=\frl_ \bbC\oplus \fru_\bbC$  is a  $\theta$-stable parabolic subalgebra of 
 $\frg_\bbC$,  with Levi factor  $\frl_\bbC$,  and if  $L$ is the normalizer in $G$ of  $\frqqq_\bbC$, 
we define the cohomological induction functor : 
\begin{align*}  \caR_{\frqqq_\bbC} &: \caM(\frg_\bbC,K)\longrightarrow  \caM(\frl_\bbC,K\cap L)     \\
X & \mapsto \Gamma^S \circ \mathrm{pro} (X \otimes \tilde \tau )  \end{align*}
where   $S =\dim(\fru_\bbC \cap \frk_\bbC)$, $\Gamma^S$ is the  $S$-th  Zuckerman derived functor  from  
$\caM(\frg_\bbC,K\cap L)$ to  $\caM(\frg_\bbC,K)$, $ \mathrm{pro}$ is the parabolic induction functor from 
 $\caM(\frl_\bbC,K\cap L)$ to  $\caM(\frg_\bbC,K\cap L)$, and  $\tilde{\tau}$ is a character of  
$L$, square root of the character  $\bigwedge^{top}(\fru_\bbC/\fru_\bbC \cap \frk_\bbC)$ (such a square root is usually 
defined only on a double cover of   $L$, but for the cases we are interested in here, {\sl i.e.} products of 
$G=\GL(n,\bbR)$, $\GL(n,\bbC)$ or  $\GL(n,\bbH)$, we can find such a square root on  $L$). 
This normalization preserves infinitesimal character.

With this notation, for   $G=\GL(2,\bbR)$, $L\simeq \bbC^\times$ and  $\fru_\bbC=\frg_\bbC^{e_1-e_2}$, we get   
\[\caR_{\frqqq_\bbC} (\gamma(x,y)) =\eta(x,y), \quad x,y \in \bbC, x-y \in \bbN.  \]
Recall the convention  $\eta(x,x)=\delta(x,0)\times \delta(x,1)$ for limits of discrete series. 
We have then  $\caR_{\frqqq_\bbC} (\gamma(x,y)) =\eta(x,y)$. 

Set  $a_0=a(\eta(x,y),n) \in M(D)$. The  standard representation  $\lambda(a_0)$ is obtained by parabolic induction 
from the representation 
\begin{small}
\[ \eta= \eta(x+\frac{n-1}{2},y+\frac{n-1}{2})\otimes  \eta(x+\frac{n-3}{2},y+\frac{n-3}{2}) \otimes \ldots \otimes 
\eta(x-\frac{n-1}{2},y-\frac{n-1}{2})\] \end{small}
of  $\GL(2,\bbR)\times \ldots \times \GL(2,\bbR)$, the representation  $\eta$ being from what has just been said 
obtained by cohomological induction from the character
\begin{small}
\[ \gamma= \gamma(x+\frac{n-1}{2},y+\frac{n-1}{2})\otimes  \gamma(x+\frac{n-3}{2},y+\frac{n-3}{2}) 
\otimes \ldots \otimes 
\gamma(x-\frac{n-1}{2},y-\frac{n-1}{2})\]\end{small}
of $\bbC^\times \times \ldots \times \bbC^\times$. Furthermore  $u(\eta(x,y),n)$ is the unique irreducible quotient 
of   $\lambda(a_0)$.

Independence of polarization results in \cite{KV}, chapter  11 show that the standard representation  $\lambda(a_0)$
 could be also obtained from  the character  $\gamma$ of  $(\bbC^\times)^n$ in the following way : 
first use parabolic induction from  $(\bbC^\times)^n$ to  $\GL(n,\bbC)$ (with respect to the usual upper triangular
 Borel subgroup) to get the standard representation 
\begin{small}\begin{equation}\label{gamma} 
 \gamma(x+\frac{n-1}{2},y+\frac{n-1}{2})\times  \gamma(x+\frac{n-3}{2},y+\frac{n-3}{2}) \times \ldots \times 
\gamma(x-\frac{n-1}{2},y-\frac{n-1}{2})\end{equation}\end{small}
 whose unique irreducible quotient is  $u(\gamma(x,y),n)$, and then the  cohomological induction functor 
$\caR_{\frqqq_\bbC}$ from  $\GL(n,\bbC)$ to  $\GL(2n,\bbR)$ (the reader can guess which  $\theta$-stable 
parabolic subalgebra  $\frqqq_\bbC$ we  use). This shows also that   $u(\delta,n)$ is the unique irreducible quotient of 
$\caR_{\frqqq_\bbC} (u(\gamma(x,y),n))$. Now, irreducibility and unitarizability theorems of 
\cite{KV} also imply, the character   $u(\gamma(x,y),n)$ of  $\GL(n,\bbC)$ being in the weakly good range, that 
$\caR_{\frqqq_\bbC} (u(\gamma(x,y),n))$ is irreducible and unitary  if  $u(\gamma(x,y),n)$ is unitary.
 Thus we get  
\[ \caR_{\frqqq_\bbC} (u(\gamma(x,y),n))= u(\eta(x,y),n) \]
and this representation is unitary if and only if  $\Re e(x+y)=0$.

In the degenerate case $x=y$ (see (\ref{limsd})), we get 
\[   \caR_{\frqqq_\bbC} (u(\gamma(x,y),n))= u(\delta(x,0),n)\times u(\delta(x,1),n).  \]

\bigskip 

\subsection{$A=\bbH$} \label{U1H}

Let  $\delta=\eta'(x,y)$, $x,y \in \bbC$, $x-y \in \bbN^*$, be an irreducible representation of  $\bbH^\times$. 
Consider the representation  $u(\eta'(x,y),n)$, and recall  the invariant  $\nu_\delta$ of definition \ref{nudelta}.
 We have seen that  $\nu_\delta=\nu$ when $x-y >1$, $\nu_\delta=\nu^2$ when  $x-y=1$. In the first case, 
the discussion for the unitarizability of  $u(\eta'(x,y),n)$ is exactly the same as in the case  $A=\bbR$: 
the standard representation  $\lambda(a_0)$ whose unique irreducible quotient is   $u(\eta'(x,y),n)$ is 
obtained by cohomological induction from $\GL(n,\bbC)$ to  $\GL(n,\bbH)$ of  the representation  $\gamma$ 
defined in  (\ref{gamma}). Furthermore    $u(\eta'(x,y),n)$ is the unique irreducible quotient of  
$R_{\frqqq'_\bbC} (u(\gamma(x,y),n))$ and is unitary if and only if  $\Re e(x+y)=0$.

When  $\nu_\delta=\nu^2$, {\sl i.e.} $x-y=1$,  we get the same results,  not for  $u(\eta'(x,y),n)$, but for 
 $\bar u(\eta'(x,y),n)$, the  Langlands quotient of the  standard representation  
\begin{small}
\begin{equation*} \eta'(x+\frac{n-1}{2},y+\frac{n-1}{2})\times  \eta'(x+\frac{n-3}{2},y+\frac{n-3}{2}) \times \ldots
 \times  \eta'(x-\frac{n-1}{2},y-\frac{n-1}{2})\end{equation*} \end{small}
\[= \nu^{\frac{n-1}{2}}\eta'(x,y)\times \nu^{\frac{n-3}{2}} \eta'(x,y) \otimes \ldots
 \times \nu^{-\frac{n-1}{2}} \eta'(x,y) .\]
Recall that    $u(\eta'(x,y),n)$ Langlands quotient of  
 \begin{align*} &\nu_\delta^{\frac{n-1}{2}}\eta'(x,y)\times \nu_\delta^{\frac{n-3}{2}} \eta'(x,y) \otimes \ldots
 \times \nu_\delta^{-\frac{n-1}{2}} \eta'(x,y)\\
&= \nu^{n-1}\eta'(x,y)\times \nu^{n-3} \eta'(x,y) \otimes \ldots
 \times \nu^{-(n-1)} \eta'(x,y). \end{align*}

With the two conditions $x-y=1$ and  $\Re e(x+y)=0$, we see that, up to a twist by a unitary character, 
we only have to study the case   $u(\eta',n)$ with $\eta'=\eta'(\frac{1}{2},-\frac{1}{2})$.
Unitarity of   $u(\eta',n)$ can be deduced from the  unitarity of the $\bar u(\eta',k)$  
 as in \cite{BR1}, using the fact that 
\begin{align}\label{ubar} \bar u (\eta',2n+1)= u(\eta',n+1)\times  u(\eta',n)\end{align}
\begin{align}\label{ubar2}\bar u (\eta',2n)= \nu^{\frac{1}{2}} u(\eta',n)\times \nu^{-\frac{1}{2}} u(\eta',n).  \end{align}

\bigskip

\section{Vogan's classification and $U(0)$ in the   archimedean case}\label{VOG}

As we have already said, $U(0)$ is established in the case $A=\bbR$ or $\bbC$ by the work of M. Baruch
filling the serious technical gap that remained in Bernstein approach. It is also possible to establish $U(0)$
from Vogan's classification, and this will work also for $A=\bbH$. Of course, this  might seem a rather convoluted
and unnatural approach, if the final goal is to prove the classification of the unitary dual in Tadi\'c's form, since
a direct comparison between the classifications is possible. But let us notice that : 

--- One of the main difficulty of Vogan's paper is to prove   some special cases of  $U(0)$ (the other difficult point is
 the  exhaustion of the list of   unitary almost spherical representations). 
 The rest of his paper uses only standard and general techniques of the representation of real reductive groups, mainly
 cohomological induction.

--- The argument which allows the comparison between the two classifications  (``independence of  polarizations'') 
is also the one leading to  $U(0)$ from Vogan's  classification.

---  There is still some hope to find an uniform proof of  $U(0)$ for all $A$.

\bigskip

In this  section, we give a brief overview of  Vogan's paper \cite{Vog2}, and how it implies  $U(0)$. Here, 
$A=\bbR,\bbC$ or $\bbH$.

Let us fix a unitary character
\[  \delta  :  \GL(1,A) \simeq  A^\times  \rightarrow \bbC^\times.   \]
It  extends canonically to a  family of unitary   characters 
\[  \delta_n  :  \GL(n,A)    \rightarrow \bbC^\times,   \]
 by composing with the determinant $ \GL(n,A)\rightarrow  \GL(1,A)$ (non commutative determinant 
of  Dieudonn\'e  if $F=\bbH$).

The   {\sl basic blocs} of Vogan's  classification are the  representations : 
\[  \nu^{i\beta} \delta_n, \quad \beta \in \bbR  \]
(with Tadi\'c's notation, $\nu^{i\beta} \delta_n=u(\nu^{i\beta}\delta,n)$ : 
 it is a unitary character of  $\GL(n,A)$), and the representations 
\[ \pi(\nu^{i\beta}\delta,n;\alpha)= \nu^{-\alpha}\nu^{i\beta}\delta_n \times  \nu^{\alpha} \nu^{i\beta}\delta_n, 
\quad 0<\alpha<\frac{1}{2}  \]
of  $\GL(2n,F)$. These are  Stein's complementary series.

Vogan considers first parabolically induced representations of the form 
\begin{equation} \label{tau} \tau = \tau_1 \times \tau_2 \times \ldots \times \tau_r \end{equation}
where each  $\tau_j$ is either a unitary character  
\[  \tau_j= \nu^{\beta_j} \delta_{n_j}, \quad \beta_j \in i\bbR, \]
or a Stein's complementary series 
\[  \tau_j= \pi(\nu^{\beta_j}\delta, n_j;\alpha), \quad  \beta_j \in i\bbR, \; 0<\alpha<\frac{1}{2}\]

The reason for these  conditions is the following  : recall our choices of maximal compact subgroups  $K(n,A)$ 
of   $\GL(n,A)$ respectively  for $A=\bbR, \bbC$ and  $\bbH$ : 
\[ \Or(n),\U(n) \text{ and } \Sp(n)  \]
and denote by  $\mu_n$ the  restriction of  $\delta_n$ to  $K(n,A)$. We say that  $\mu_n$ is a 
 {\sl special} $1$-dimensional   representation of  $K(n,A)$.  If  $A=\bbR$, since  $\mu_n$  factorizes through the  
 determinant, there are two special  representations of  $\Or(n)$ : the trivial  representation, 
and the  sign of the  determinant. If  $A=\bbC$, special  representations of  $\U(n)$ are obtained by composing 
the determinant (with  values in  $\U(1)$), and a character of  $\U(1)$ (given by an integer). Finally, if  
 $A=\bbH$ the only  special representation of  $\Sp(n)$ is the trivial one.
A  representation of   $\GL(n,A)$ is said to be  {\sl almost spherical} (of  type $\mu_n)$ 
 if it contains the special   $K$-type  $\mu_n$.
This generalizes  spherical representations. The characters  $\delta_{n} \nu^\beta$ are  exactly  the ones  whose 
 restriction to $K(n,A)$ is  $\mu_n$.  The $\tau_i$'s  above are thus either almost spherical  unitary 
characters of type  $\mu$  (the family $\mu=(\mu_n)_n$ is fixed), or  almost  spherical Stein's complementary series 
 of type $\mu$. 

Then Vogan shows  the following  (\cite{Vog2}, Theorem 3.8): 

\begin{thm}
The   representations  $ \tau = \tau_1 \times \tau_2 \times \ldots \times \tau_r$ are 

$(i)$ unitary

$(ii)$ irreducible

Furthermore, every irreducible, almost spherical of type $\mu$, unitary  representation is obtained in this way, and  
two  irreducible, almost spherical of type $\mu$, unitary representations   
\[ \tau = \tau_1 \times \tau_2 \times \ldots \times \tau_r \]
and 
 \[ \tau' = \tau_1' \times \tau_2' \times \ldots \times \tau_s' \]
are  equivalent if and only if  the  multisets $\{\tau_i'\}$ and   $\{\tau_j\}$ are equal. 
\end{thm}  

Let us notice that this theorem contains a special case of $U(0)$ :  this is the point $(ii)$. It can be proved
using  Proposition 2.13 in \cite{Ba1} and results of S. Sahi (\cite{Sa}, Thm 3A).

Furthermore, the  classification of irreducible, almost spherical, unitary  
representations it gives  coincide  with  Tadi\'c's classification. (One has to notice that 
an  irreducible, almost spherical, unitary   representation is such with respect to an unique 
special  $K$-type : special $K$-types are minimal, and minimal   $K$-types for  $\GL(n,A)$ are unique,
 and appear with 
 multiplicity $1$).

 Vogan classification of the unitary dual of  $\GL(n,A)$ reduces matters to this particular case of 
 almost spherical  representations using cohomological induction functors preserving irreducibility and unitarity. 
More precisely, let us recall some material about Vogan's classification of the admissible dual 
of a real reductive group by $G$ by minimal $K$-types   (\cite{Vog1}). 
To each irreducible  representation of  $G$ is attached a finite number of minimal  $K$-types. 
As we said above, for  $G=\GL(n,A)$, the minimal $K$-type is unique, and appears with  multiplicity $1$. 
This gives a partition  (which can be explicitly given in terms of  Langlands classification) of the 
admissible dual of $\GL(n,A)$. 

Vogan's classification of the unitary dual deals  with each term of this partition separately. 
To each  irreducible representation $\mu$ of the  compact group  $K(n,A)$ is attached a subgroup
 $L$ of  $\GL(n,A)$ with maximal compact subgroup  $K_L:=K(n,A)\cap L$, and 
an  irreducible representation $\mu_L$  of  $K_L$.  The subgroup  $L$ is a product of groups of the form
 $\GL(n_i,A_i)$,  \[ K(n,A)\cap L \simeq \prod_i K(n_i,A_i)\   \]
 and  $\mu_L$ is a tensor product of special representations of the   $K(n_i,A_i)$. 

As opposed to Tadi\'c's classification  which uses only parabolic induction  functors, 
Vogan's classification of   $\GL(n,\bbR)$ for instance, will use classification of the almost spherical unitary 
dual   of groups   $\GL(k,\bbC)$. More precisely :

- For  $F=\bbR$, the subgroups  $L$ are products of $\GL(k,\bbR)$ and  $\GL(m,\bbC)$. 

- For  $F=\bbC$, the subgroups $L$ are products of  $\GL(k,\bbC)$. 

- For  $F=\bbH$,  the subgroups  $L$ are products of $\GL(k,\bbH)$ and  $\GL(m,\bbC)$. 

A combination of parabolic and cohomological induction functors then defines a functor 
 \[ \caI_L^G     \]
from  $\caM(L)$ to  $\caM(\GL(n,A))$ with the following properties : 

- $\caI_L^G$ sends an irreducible (resp. unitary) representation of  $L$ with minimal  $K_L$-type $\mu_L$
to an irreducible (resp. unitary) representation of  $\GL(n,F)$ with minimal  $K$-type  $\mu$.

 - $\caI_L^G$ realizes  a bijection between equivalence classes of irreducible  unitary representations of  $L$
with minimal  $K_L$-type $\mu_L$ and equivalence classes of irreducible  unitary representations of  $\GL(n,F)$. 
with minimal  $K$-type $\mu$.

\bigskip 
 
>From this point of view,  to establish  $U(0)$, the first thing to do is to check that 
products of representations of the form  (\ref{tau}) for different families of special $K$-types  $\mu$ are  
 irreducible. For  $F=\bbH$, there is nothing to check since there is only one  family of special  
$K$-types  $\mu=(\mu_n)_n$. For  $F=\bbR$, there are two  
families of special   $K$-types, the trivial and sign characters of the  determinant of  $\Or(n)$. 
The relevant result is then lemma 16.1 of  \cite{Vog2}. For  $F=\bbR$, we have now obtained all 
irreducible unitary representations  which are products of  $u(\delta,k)$ and  $\pi(\delta,k;\alpha)$
with $\delta$ any unitary character of  $\GL(1,\bbR)=\bbR^\times$. 

The  case $A=\bbC$ is simpler and dealt with as follows. 
Let us notice first that since square integrable modulo center representations of  
$\GL(n,\bbC)$ exist only for $n=1$, the above assertion shows that  we get all representations of  Tadi\'c's
 classification, and establishes  $U(0)$. In that case, the subgroups  $L$ from which we use
cohomological induction are of the form
\[ L=\GL(n_1,\bbC)\times \ldots \GL(n_r,\bbC). \]
The cohomological induction setting is that  $\frl_\bbC=\mathrm{Lie}(L)_\bbC$ 
is a Levi factor of a $\theta$-stable parabolic subalgebra  $\frqqq_\bbC$ 
of  $\frg_\bbC= \mathrm{Lie}(\GL(n,\bbC))_\bbC$. But  $L$ is also a Levi factor 
of a parabolic subgroup of  $\GL(n,\bbC)$. Thus there are two ways of inducing from $L$ to 
 $\GL(n,\bbC)$: parabolic and cohomological induction. An 
`independence of  polarization' result  (\cite{Vog2}, Theorem 17.6, see  \cite{KV}, Chapter  11 for a proof), 
asserts that the two coincide. This finishes the case $A=\bbC$. 

Let us now finish to discuss the cases  $A=\bbR$ and   $A=\bbH$.  Representations  
from Tadi\'c's classification which are still missing are the ones built from  $u(\delta,k)$'s and 
 $\pi(\delta,k;\alpha)$'s  with $\delta$  a square integrable modulo center representation of 
 $\GL(2,\bbR)$ or  $\bbH^\times$. As we have seen in  \ref{U1R}  a square integrable modulo center representation of 
  $\GL(2,\bbR)$ or $\bbH^\times$ is obtained by cohomological induction from  the subgroup 
 $L \simeq \bbC^\times$ of   $\GL(2,\bbR)$ or $\GL(1,\bbH)=\bbH^\times$.
This explains somehow why cohomological induction will produce the missing representations.
Let us explain this  : 

--- case  $F=\bbR$ : we start with representations of the form
\[ u(\chi_a,k_a), \; \pi(\chi_b,k_b;\alpha_b),\; u(\chi_c,k_c), \; \pi(\chi_d,k_d;\alpha_d), 
\; u(\chi_e,k_e), \; \pi(\chi_f,k_f;\alpha_f))   \]
where  $u(\chi_a,k_a)$ are unitary characters of  $\GL(k_a,\bbC)$,   $\pi(\chi_b,k_b;\alpha_b)$ 
are Stein complementary series of  $\GL(2k_b,\bbC)$,  $u(\chi_c,k_c)$ are unitary characters of 
$\GL(k_c,\bbR)$ of trivial type  $\mu$,   $\pi(\chi_d,k_d;\alpha_d)$ 
are Stein complementary series of  $\GL(2k_d,\bbR)$ of trivial type $\mu$, 
 $u(\chi_e,k_e)$  are unitary characters of $\GL(k_c,\bbR)$ of  type $\mu=\mathrm{sgn}$, 
 $\pi(\chi_f,k_f;\alpha_f)$ are Stein complementary series of  $\GL(2k_f,\bbR)$ of type $\mu=\mathrm{sgn}$.

The tensor product
\begin{small}
\[ \bigotimes_a u(\chi_a,k_a)\bigotimes_b\pi(\chi_b,k_b;\alpha_b)\bigotimes_c  u(\chi_c,k_c) \bigotimes_d 
\pi(\chi_d,k_d;\alpha_d) \bigotimes_e u(\chi_e,k_e) \bigotimes_f \pi(\chi_f,k_f;\alpha_f))  \]\end{small}
is a representation of the Levi subgroup
\begin{small}\[ \prod_a \GL(k_a,\bbC) \prod_b \GL(2k_b,\bbC)\prod_c \GL(k_c,\bbR) \prod_d 
\GL(2k_d,\bbR)\prod_e \GL(k_e,\bbR)
 \prod_f \GL(2k_f,\bbR) \]\end{small}
of  $\GL(n,\bbR)$, where $n=\sum_a 2k_a+\sum_b 4k_b+\sum_c k_c +\sum_d 2k_d+ \sum_e k_e +\sum_k 2k_f$.  

As we saw, we first form almost spherical representations of a given type by parabolic induction. Thus we induce
\[\bigotimes_c  u(\chi_c,k_c) \bigotimes_d  \pi(\chi_d,k_d;\alpha_d)\]
from 
\[ \prod_c \GL(k_c,\bbR) \prod_d \GL(2k_d,\bbR)\]
to  $\GL(q_0,\bbR)$, where  $q_0=\sum_c k_c +\sum_d 2k_d$, 
obtaining an irreducible unitary  spherical  representation   $\pi_0$, and similarly   
\[\bigotimes_e  u(\chi_e,k_e) \bigotimes_f  \pi(\chi_f,k_f;\alpha_f)\]
from 
\[ \prod_e \GL(k_e,\bbR) \prod_d \GL(2k_f,\bbR)\]
to  $\GL(q_1,\bbR)$, where  $q_1=\sum_e k_e +\sum_f 2k_f$, obtaining an irreducible unitary almost 
 spherical  type  $\mu=\mathrm{sgn}$  representation. 

Then we mix spherical and almost spherical of  type  $\mu=\mathrm{sgn}$ representations inducing   
parabolically  $\pi_0\times \pi_1$ from   $\GL(q_0,\bbR)\times \GL(q_1,\bbR)$ to  $\GL(q_0+q_1,\bbR)$ : we get 
  an irreducible unitary    representation  $\pi$ of  $\GL(q_0+q_1,\bbR)$. 

The  group  $\prod_a \GL(k_a,\bbC) \prod_b \GL(2k_b,\bbC) \times \GL(q_0+q_1,\bbR)$ is denoted by  $L_\theta$ in 
\cite{Vog2}. 
Applying cohomological induction functor  $\caI_{L_\theta}^G$ to the  representation 
\[  \bigotimes_a u(\chi_a,k_a)\bigotimes_b\pi(\chi_b,k_b;\alpha_b)\otimes \pi  \]
of  $L_\theta$, we get  an irreducible unitary    representation    $\rho$ of  $\GL(n, \bbR)$.

Independence of  polarization theorems  \cite{Vog2}, Theorem 17.6, Theorem 17.7 and 17.9 
(see \cite{KV}, Chapter  11),
allows us to invert the order of the two types of induction. We could in fact start with  cohomological
 induction, inducing each   \[ u(\chi_a,k_a) \]
from  $\GL(k_a,\bbC)$ to  $\GL(2k_a,\bbR)$. 
In non degenerate case, following the terminology of 
\cite{Vog2}, definition 17.3,  we get representations  $u(\delta_a,2k_a)$, where  
$\delta_a$ is a square integrable modulo center irreducible representation of  $\GL(2,\bbR)$. 
In the degenerate case, $\delta_a$ 
is a limit of discrete series (\ref{limsd}). These are almost spherical representations that we had before
 (see \cite{Vog2}, prop. 17.10).

In the same way, we induce all 
\[ \pi(\chi_b,k_b;\alpha_b)\]
from   $\GL(2k_b,\bbC)$ to  $\GL(4k_b,\bbR)$. In the  non degenerate case, 
  we get representations $\pi(\delta_b,2k_b;\alpha_b)$, where  
$\delta_b$ is as above. In the degenerate case, we still get almost spherical representations.

 The parabolically induced representation from 
\[ \prod_a \GL(2k_a,\bbR) \prod_b \GL(4k_b,\bbR) \times \GL(q_0+q_1,\bbR) \]
to  $\GL(n,\bbR)$ of  
\[  \bigotimes_a u(\delta_a,k_a)\bigotimes_b\pi(\delta_b,k_b;\alpha_b)\otimes \pi  \]
is  $\rho$ (and thus  irreducible), see \cite{Vog2}, Theorem 17.6 .

This finishes the comparison of the two  classifications. The case  $A=\bbH$ is entirely similar. 

We deduce  $U(0)$  using again independence of polarization. 
We want to show that  $\rho=\rho_1 \times \rho_2$ is irreducible if  $\rho_1$ and  $\rho_2$ are  
irreducible and  unitary. We write  $\rho_1$ and  $\rho_2$ as above using first  
 cohomological induction and then,  parabolic induction. Using parabolic induction by stage, we see that 
 $\rho_1 \times \rho_2$ is also written in this form. Using again  independence of polarization
 we write $\rho$ as a parabolically then cohomologically induced representation,
and we see that as such, this is a representation appearing in Vogan's classification which is therefore irreducible.

\section{ Jacquet-Langlands correspondence in the archimedean case}\label{JaLa}

Ideas in this section are taken from \cite{AH} which deals with a similar problem (Kazhdan-Patterson lifting).

\subsection{Jacquet-Langlands and coherent families}

Since we need to consider simultaneously the case $A=\bbR$ and $A=\bbH$, we 
add  relevant superscripts to the notation when needed as in Section \ref{JL1}. 
We have noticed  that Jacquet-Langlands correspondence between essentially square integrable
modulo center irreducible representations of $\GL(2,\bbR)$ and irreducible representations of  
$\bbH^\times$ is given at the level of  Grothendieck groups  by 
\[ \JL(\eta(x,y))= - \eta'(x,y) \]
Representations in $D_1$ are sent to $0$. 
We extend this linearly to an algebra morphism :
\[ \caR^\bbR  \rightarrow \caR^\bbH. \]

\begin{lemme}
 Jacquet-Langlands correspondence preserves  supports.
\end{lemme}
\pr  $a \in M(D)$, $a=(\eta(x_1,y_1), \ldots ,\eta(x_r,y_r))$. We have then  
\[ \JL(\lambda(a))= (-1)^r \lambda(a') \]
where  $a=(\eta'(x_1,y_1), \ldots ,\eta'(x_r,y_r))$.
The  support of $a$ is  $(x_1,y_1,\ldots ,x_r,y_r)$, and this is also the support of  $a'$. 
\bigskip

We recall now the definition of a coherent family of  Harish-Chandra modules.  

\begin{defi}
Let  $G$ be a real reductive group, $H$  a  Cartan subgroup,  $\frg_\bbC$ and $\frh_\bbC$  the respective 
 complexification
 of their Lie algebras and   $\Lambda$ the  lattice of weights of  $H$ in finite dimensional representations of  $G$.
 A coherent family of (virtual)  Harish-Chandra modules   based at  $\lambda \in \frh_\bbC^*$
is a family  $$\{ \pi(\lambda+\mu)\, |\,  \mu \in \Lambda  \}$$ 
($\lambda+\mu$ is just a formal symbol, since the two terms are not in the same group)
in the Grothendieck group $\caR(G)$ such that

---  The infinitesimal character of  $ \pi(\lambda+\mu)$ is given by  $\lambda+d \mu$.
 
--- For all finite dimensional representation  $F$ of $G$, we have, with  $\Delta(F)$ denoting the set of weights of 
$H$ in  $F$, the following identity in   $\caR(G)$ : 
 \[  \pi(\lambda+\mu)\otimes F = \sum_{\gamma \in \Delta(F)} \pi(\lambda+\mu+\gamma). \]
\end{defi}

Jacquet-Langlands correspondence preserves coherent families :

\begin{lemme}Let us  identify two Cartan subgroups  $H$ and  $H'$ respectively of 
  $\GL(2n,\bbR)$ and   $\GL(n,\bbH)$ isomorphic to $(\bbC^\times)^n$.
Let  $\pi(\lambda+\mu)$ be a coherent family of  Harish-Chandra modules for  $\GL(2n,\bbR)$ based 
at $\lambda \in \frh_\bbC^*$.
Then  $\JL(\pi(\lambda+\mu))$ is a coherent family for    $\GL(n,\bbH)$.
\end{lemme}
\pr The first property of coherent families is satisfied by  $\JL(\pi(\lambda+\mu))$ because of the previous lemma.
For the second property, let us remark first that 
 $\GL(2n,\bbR)$ and  $\GL(n,\bbH)$ being two real forms of 
$\GL(2n,\bbC)$, a finite dimensional representation  $F$ of one of these two groups is in fact  the restriction of 
a finite dimensional representation of  $\GL(2n,\bbC)$.
 We get for all regular element  $g'$ of  $\GL(n,\bbH)$  corresponding to an 
element  $g$ in  $\GL(2n,\bbR)$,  
\begin{align*}  &\sum_{\gamma \in \Delta(F)} \Theta_{\JL( \pi(\lambda+\mu+\gamma))}(g') =
 \sum_{\gamma \in \Delta(F)} \Theta_{\pi(\lambda+\mu+\gamma)}(g)
=  \Theta_{\pi(\lambda+\mu)\otimes F}(g)\\
&=  \Theta_{\pi(\lambda+\mu)}(g) \, \Theta_F(g)
=\Theta_{\JL(\pi(\lambda+\mu))}(g') \Theta_F(g')
=\Theta_{\JL(\pi(\lambda+\mu))\otimes F}(g'), 
\end{align*}
so  $\sum_{\gamma \in \Delta(F)} \JL(\pi(\lambda+\mu+\gamma))=\JL(\pi(\lambda+\mu))\otimes F$. 
\qed 

\bigskip

\subsection{Jacquet-Langlands and cohomological induction}

The cohomological induction functor  $\caR_{ \frqqq_\bbC}$ introduced in  \ref{U1R}
preserves  irreducibility and  unitarity when the infinitesimal character of the 
induced module  satisfies certain positivity properties 
with respect to  $ \frqqq_\bbC$ (``weakly good range''). Furthermore with the same conditions, other derived functors 
 $\Gamma^i(\mathrm{pro}(\bullet \otimes \tilde \tau))$, $i\neq S$ vanish.
This is not true in general, and this is the reason why we need to consider Euler-Poincar\'e characteristic : 
 \[  \widehat \caR_{ \frqqq_\bbC}:= \sum_{i} (-1)^i \Gamma^i(\mathrm{pro}(\bullet \otimes \tilde \tau)).  \]
This is not a functor between  $\caM(L)$ and  $\caM(G)$ anymore, but simply a morphism between  the Grothendieck 
groups $\caR(L)$ and $\caR(G)$.

\begin{lemme}
The  morphism $ \widehat \caR_{ \frqqq_\bbC}: \; \caR(L) \rightarrow \caR(G)$ preserves  coherent families.
\end{lemme}

\pr  The functors   $\Gamma^i(\mathrm{pro}(\bullet \otimes \tilde \tau))$ are normalized in 
order to preserve  infinitesimal character, and thus the first property of coherent family is preserved.

Let  $\pi(\lambda+\mu)$ be a coherent family of Harish-Chandra for  $(\frl,L\cap K)$.
We want to show that for any finite dimensional representation  $F$  of $G$, 
\begin{equation}
 \widehat \caR_{ \frqqq_\bbC}(\pi(\lambda+\mu)) \otimes F=
 \sum_{\gamma \in \Delta(F)}  \widehat \caR_{ \frqqq_\bbC}(\pi(\lambda+\mu+\gamma))
\end{equation}

But  \begin{align*}
& \sum_{\gamma \in \Delta(F)} \widehat  \caR_{ \frqqq_\bbC}(\pi(\lambda+\mu+\gamma))= \widehat   \caR_{ \frqqq_\bbC}
\left(\sum_{\gamma \in \Delta(F)} \pi(\lambda+\mu+\gamma) \right)\\
=&  \widehat   \caR_{ \frqqq_\bbC} (\pi(\lambda+\mu) \otimes F)
\end{align*}

It is then enough to show that for any   $(\frl,L\cap K)$-module $X$, 
 \begin{equation} \label{OF} \widehat   \caR_{ \frqqq_\bbC} (X)\otimes F =  \widehat   \caR_{ \frqqq_\bbC} (X\otimes F ) 
\end{equation}

Let  $U$ be a  any  $(\frg, K)$-module. Let us compute, using adjunction properties of the functors involved: 
 \begin{align*}
&\Hom_{\frg,K} (U,\Gamma( \mathrm{pro}((X\otimes F) \otimes \tilde \tau)))\simeq 
 \Hom_{\frl,L\cap K} (U,X\otimes F \otimes \tilde \tau)\\
&\simeq \Hom_{\frl,L\cap K} (U,X\otimes (F^*)^* \otimes \tilde \tau)
\simeq \Hom_{\frl,L\cap K} (U, \Hom_\bbC (F^*,X\otimes  \tilde \tau)\\
&\simeq \Hom_{\frl,L\cap K} (U \otimes F^*, X\otimes  \tilde \tau))
 \simeq \Hom_{\frg,K} (U \otimes F^*,\Gamma( \mathrm{pro}(X\otimes \tilde \tau)))\\
& \simeq \Hom_{\frg,K} (U,\Gamma( \mathrm{pro}(X\otimes \tilde \tau))\otimes F)
\end{align*}
We deduce from this that  $\Gamma( \mathrm{pro}(X\otimes \tilde \tau \otimes F))\simeq 
\Gamma( \mathrm{pro}(X\otimes \tilde \tau))\otimes F$.

The same is true for $\Gamma^i$ replacing $\Gamma$ in the computation above. This can be seen 
using general arguments using  the exactness of the functor  $\bullet \otimes F$.
Thus,  for all $i\geq 0$, 
$\Gamma^i( \mathrm{pro}(X\otimes   \otimes \tilde\tau \otimes F))
\simeq \Gamma^i( \mathrm{pro}(X\otimes \tilde \tau))\otimes F$,
 which implies 
(\ref{OF}).

\qed
\bigskip

Let us now denote   $\widehat \caR_{ \frqqq_\bbC}^\bbR$ and  $\widehat \caR_{ \frqqq'_\bbC}^\bbH$ the 
Euler-Poincar\'e  morphisms  of  cohomological induction  between  $\GL(1,\bbC)$ and respectively 
 $\GL(2, \bbR)$ and  $\GL(1,\bbH)$, where 
 $\frqqq_\bbC$ and  $\frqqq_\bbC'$ are as  \ref{U1R} and   \ref{U1H}.

\begin{lemme}\label{Cohcont}
With the notation above, and  $x,y \in \bbC$, $x-y \in \bbZ$,
\[ \JL(\widehat \caR_{ \frqqq_\bbC}^\bbR(\gamma(x,y))) =- \widehat \caR_{ \frqqq'_\bbC}^\bbH(\gamma(x,y)) \] 
\end{lemme}
\pr When  $x-y\geq 0$, we have  $$\widehat \caR_{ \frqqq_\bbC}^\bbR(\gamma(x,y))=- \caR_{ \frqqq_\bbC}^\bbR(\gamma(x,y))=
-\eta(x,y)$$
 and 
$$\widehat \caR_{ \frqqq'_\bbC}^\bbH(\gamma(x,y))= - \caR_{ \frqqq'_\bbC}^\bbH(\gamma(x,y))=-\eta'(x,y).$$
The formula is thus true in this case. The case $x-y<0$ follows because  
$\JL(\widehat \caR_{ \frqqq_\bbC}^\bbR(\gamma(x-n,y+n)))$ and   $\widehat \caR_{ \frqqq'_\bbC}^\bbH(\gamma(x-n,y+n))$ 
are two coherent families which coincide for  $n \geq 0$, and are therefore equal.

\begin{thm}
Let  $\caR_{ \frqqq_\bbC}^\bbR$ and  $\caR_{ \frqqq'_\bbC}^\bbH$ be the  cohomological induction functors 
from  $\GL(n,\bbC)$ to respectively  $\GL(2n,\bbR)$ and   $\GL(n,\bbH)$.
We have then  \[\JL\circ  \widehat \caR_{ \frqqq_\bbC}^\bbR= (-1)^n \widehat \caR_{ \frqqq'_\bbC}^\bbH .\]
\end{thm}

\pr 
It is enough to show that the formula holds on the basis  $\lambda(a)$, $a \in M(D)$ of  $\caR^\bbC)$.
Let  $a \in M(D)$, $a=(\gamma(x_1,y_1), \ldots ,\gamma(x_r,y_r))$. We compute
\begin{align*}
\JL\circ  \widehat \caR_{ \frqqq_\bbC}^\bbR(\lambda(a))&=\JL\circ  \widehat \caR_{ \frqqq_\bbC}^\bbR(
\gamma(x_1,y_1)\times  \ldots \times \gamma(x_r,y_r))\\
&=\JL(i_{\GL(2,\bbR)^r}^{\GL(2r,\bbR)} \circ \widehat \caR_{ \frqqq_\bbC}^\bbR(
\gamma(x_1,y_1)\otimes  \ldots \otimes \gamma(x_r,y_r))\\
&=i_{\GL(1,\bbH)^r}^{\GL(r,\bbH)}\circ  \JL(\widehat \caR_{ \frqqq_\bbC}^\bbR(
\gamma(x_1,y_1)\otimes  \ldots \otimes \gamma(x_r,y_r))\\
&= (-1)^r i_{\GL(1,\bbH)^r}^{\GL(r,\bbH)}\circ \widehat \caR_{ \frqqq'_\bbC}^\bbH (
\gamma(x_1,y_1)\otimes  \ldots \otimes \gamma(x_r,y_r))\\
&= (-1)^r \widehat \caR_{ \frqqq'_\bbC}^\bbH (
\gamma(x_1,y_1)\times  \ldots \times \gamma(x_r,y_r))\\
&= (-1)^r \widehat \caR_{ \frqqq'_\bbC}^\bbH (\lambda(a))
\end{align*}
We have used independence of polarization theorem of \cite{KV},  to replace a part of  cohomological induction  
by parabolic  induction, and the fact that  $\JL$ commutes with parabolic induction. \qed

\begin{cor} Recall the representations $\bar u(\eta',n)$ introduced in  \ref{U1H}.
We have  $$\JL(u(\eta(x,y),n)= (-1)^n \; \bar u(\eta'(x,y),n),$$
 $x,y\in \bbC$, $x-y \in \bbN$. 
\end{cor}

Recall that when $x-y\neq 1$, then $\bar u(\eta'(x,y),n)= u(\eta'(x,y),n)$ (see \ref{U1H}).

\pr This follows from the theorem and the formulas   $\caR_{ \frqqq_\bbC}^\bbR(u(\gamma(x,y))=u(\eta(x,y),n)$, 
  $\caR_{ \frqqq'_\bbC}^\bbH(u(\gamma(x,y))=\bar u(\eta'(x,y),n)$ obtained in 
  \ref{U1R} and  \ref{U1H}. \qed 

\bigskip

To be able to compute the transfer to  $\GL(n,\bbH)$ of any irreducible unitary representation of  
$\GL(2n,\bbR)$, we need to compute the transfer of the  $u(\delta,k)$ when  $\delta \in D_1^\bbR$.
But, in this case, if  $\delta=\delta(\alpha,\epsilon)$,
\[ u(\delta(\alpha,\epsilon),2k) = \delta(\alpha,\epsilon) \circ  \det,   \]
and we know from  \cite{DKV}  that the  transfer of this  character is the character 
\[  \delta(\alpha,\epsilon) \circ  RN  \]
($RN$ is the reduced  norm) which is  
\[  u(\eta'(\alpha+\frac{1}{2}, \alpha-\frac{1}{2}),k).  \]

 From this, we get 
\begin{thm}\label{archimtransfer}
Let  $u$ be an irreducible unitary representation of  $\GL(2n,\bbR)$. Then  
$\JL(u)$ is either $0$, or up to a sign, an irreducible unitary representation of  $\GL(n,\bbH)$.
For representations  $u(\delta,k)$, we get:

--- if  $\delta=\delta(\alpha,\epsilon) \in D_1^\bbR$, 
\[  \JL(u(\delta(\alpha,\epsilon),2k)) =   u(\eta'(\alpha+\frac{1}{2}, \alpha-\frac{1}{2}),k)   \]
 
--- if  $\delta=\eta(x,y) \in D_2^\bbR$, \[ \JL(u(\eta(x,y)),k)= (-1)^k  \bar  u(\eta'(x,y),k). \]
\end{thm}

To make it simple, a character is sent by $\JL$ on the corresponding character, while if 
$\delta\in D_2^F$ and $\delta'=\mathbf{C}(\delta)=-\JL(\delta)$,  
 then $\JL(u(\delta, k))=(-1)^k \bar{u}(\delta',k)$.

 In the first case note that we deal with a slightly different situation from non archimedean fields, 
since the reduced norm of $\bbH$ is {\it not} surjective, but has image in $\bbR_+^*$. In particular, 
if $s$ is the character sign of the determinant on $GL_{2k}(\bbR)$, then $\JL(s)$ is the trivial 
character of $GL_k(\bbH)$. In the non archimedean case, it is easy to check that $\JL$ is injective 
on the set of representations $u(\delta,k)$. 

The above theorem gives a correspondence   between irreducible unitary representations of $\GL(2n,\bbR)$
and of  $\GL(n,\bbH)$, by forgetting the signs. As in the introduction, we denote this correspondence  by 
$|\LJ|$. Using (\ref{ubar}) and (\ref{ubar2}), we easily reformulate
 the result as in the introduction.

\section{Character formulas and ends  of complementary series}\label{ends}

From Tadi\'c's classification of the unitary dual, and the character formula for induced representations, 
the character of any irreducible unitary representation of  $\GL(n,A)$ can be computed from the characters of the 
 $u(\delta,n)$, $\delta \in D$, $n \in \bbN$. It is remarkable that the characters of the $u(\delta,n)$ 
can be computed, or more precisely, expressed in terms of characters of square integrable modulo center
 representations. We give also composition series of ends of complementary series. This information is important 
for the topology of the unitary dual (see \cite{Ta6}).

\subsection{$A=\bbC$}

Let $\gamma= \gamma(x,y)$ be a character of  $\bbC^\times$, $x,y \in \bbC$, $x-y=r\in \bbZ$. 
The  representation $u(\gamma(x,y),n)$ is the character 
\[ \det \circ \gamma \]
of  $\GL(n,\bbC)$.
There is a formula, due to Zuckerman, for the trivial character of any real reductive group, obtained from 
a finite length resolution of the trivial representation by standard modules in the category $\caM(G)$.

For  $\GL(n,\bbC)$, this formula is, denoting  $\mathbf{1}_{\GL(n,\bbC)}$ the trivial representation
 
\begin{equation*}
\mathbf{1}_{\GL(n,\bbC)}=u(\gamma(0,0),n)= \sum_{w \in \frS_n} (-1)^{l(w)} \prod_{i=1}^n
\gamma(\frac{n-1}{2}-i+1,\frac{n-1}{2}-w(i)+1 ) 
\end{equation*}

From this, we get by tensoring with  $\gamma(x,y)$, 
\begin{equation}\label{Zu}
u(\gamma(x,y),n)= \sum_{w \in \frS_n} (-1)^{l(w)} \prod_{i=1}^n
\gamma(x+\frac{n-1}{2}-i+1,y+\frac{n-1}{2}-w(i)+1 ) 
\end{equation}

Set  $\gamma_{i,j}= \gamma(x+\frac{n-1}{2}-i+1,y+\frac{n-1}{2}-j+1 ) \in \caR$. The formula above becomes :
\begin{equation}\label{Zu2}
u(\gamma(x,y),n)= \det((\gamma_{i,j})_{1\leq i,j \leq n} )
\end{equation}

From the Lewis Carroll identity (\cite{CR}), we deduce easily from this a formula for composition series  
of ends of complementary series. This was obtained previously by Tadi\'c \cite{Ta3}, using partial results of
 Sahi \cite{Sah},  but the proof was complicated.
For an easy formula, set 
\[ \gamma(x,y)=\delta(\beta,r)\]
with, $r=x-y$, $2\beta=x+y$

\begin{prop} With the notation above, and  $n \geq 2$
\begin{align}\label{boutC}
& \nu^{-\frac{1}{2}} u(\delta(\beta,r) ,n)\times \nu^{\frac{1}{2}} u(\delta(\beta,r),n) \\
\nonumber&  \quad \quad =  u(\delta(\beta,r),n+1)\times  u(\delta(\beta,r),n-1)\\
\nonumber &  \quad \quad \quad \quad \quad +  u(\delta(\beta,r+1),n)\times  u(\delta(\beta,r-1),n)
\end{align}
\end{prop}

\subsection{$A=\bbR$}

Let $\eta(x,y)$ be an essentially square integrable modulo center representation  of $\GL(2,\bbR)$, 
$x,y \in \bbC$, $x-y=r\in \bbN^*$. Since 
\[ u(\eta(x,y),n)=  - \widehat \caR_{\frqqq_\bbC}^\bbR( u(\gamma(x,y))),  \]
we get from (\ref{Zu}) that 
\[  u(\eta(x,y),n)=  - \sum_{w \in \frS_n} (-1)^{l(w)} \prod_{i=1}^n
\widehat \caR_{\frqqq_\bbC}^\bbR(\gamma(x+\frac{n-1}{2}-i+1,y+\frac{n-1}{2}-w(i)+1 )).  \]

We have noticed in the proof of Lemma \ref{Cohcont} that 
$-\widehat \caR_{\frqqq_\bbC}^\bbR(\gamma(x-n,y+n))$ is a coherent family of representation of 
$\GL(2,\bbR)$ such that $-\widehat \caR_{\frqqq_\bbC}^\bbR(\gamma(x-n,y+n))=\eta(x-n,y+n)$ when $x-n > y+n$. Set  
 $\tilde \eta (x-n,y+n)=-\widehat \caR_{\frqqq_\bbC}^\bbR(\gamma(x-n,y+n))$. Then we get 
\[  u(\eta(x,y),n)=  (-1)^{n+1} \sum_{w \in \frS_n} (-1)^{l(w)} \prod_{i=1}^n \tilde 
\eta(x+\frac{n-1}{2}-i+1,y+\frac{n-1}{2}-w(i)+1 ).  \]

Set  $\tilde \eta_{i,j}= \tilde \eta (x+\frac{n-1}{2}-i+1,y+\frac{n-1}{2}-j+1 )$. The formula above becomes :
\begin{equation}\label{Zu3}
u(\eta(x,y),n)= (-1)^{n+1} \det((\tilde \eta _{i,j})_{1\leq i,j \leq n} )
\end{equation}
Again  from the Lewis Carroll identity (\cite{CR}), we deduce easily from this a formula for composition series  
of ends of complementary series

\begin{prop} With the notation above, $n \geq 2$,  $x-y>1$,
\begin{align}\label{boutR}
& \nu^{-\frac{1}{2}} u(\eta(x,y) ,n)\times \nu^{\frac{1}{2}} u(\eta(x,y),n) \\
\nonumber &  \quad \quad =  u(\eta(x,y),n+1)\times  u(\eta(x,y),n-1)\\
\nonumber &  \quad \quad \quad \quad \quad +  u(\eta(x+\frac{1}{2},y-\frac{1}{2}),n)\times  
 u(\eta(x-\frac{1}{2},y+\frac{1}{2}),n).
\end{align}

If $x=y+1$, recall the convention that 
\[ \eta(x-\frac{1}{2},x-\frac{1}{2})=\delta(x-\frac{1}{2}, 0)\times \delta(x-\frac{1}{2},1).\]
We get 
\begin{align}\label{boutR1}
& \nu^{-\frac{1}{2}} u(\eta(x,x-1) ,n)\times \nu^{\frac{1}{2}} u(\eta(x,x-1),n) \\
\nonumber&  \quad \quad =  u(\delta(x,x-1),n+1)\times  u(\eta(x,x-1),n-1)\\
\nonumber&  \quad \quad \quad \quad \quad  + u(\eta(x+\frac{1}{2},x-\frac{3}{2}),n)\times  
[ u(\delta(x-\frac{1}{2},0),n)\times  
 u(\delta(x-\frac{1}{2},1),n)]
\end{align}

\end{prop}

\begin{rmq} We cannot deduce by our method the  composition series of the ends of complementary series
for $u(\delta,n)$ when $\delta \in D_1$. There is still a formula for the character of $u(\delta,n)$, 
since  $u(\delta,n)=\delta \circ \det$ is a one-dimensional representation (Zuckerman), but no interpretation
for the right-hand-side of  this formula as a determinant, so we cannot apply the Lewis Carroll identity.  
\end{rmq}

\subsection{$A=\bbH$} The discussion is similar to the real case for the $u(\eta'(x,y),n)$ when $x-y\geq 2$.

\begin{prop} With the notation above, $n \geq 2$,  $x-y\geq 2$,
\begin{align}\label{boutH}
& \nu^{-\frac{1}{2}} u(\eta'(x,y) ,n)\times \nu^{\frac{1}{2}} u(\eta'(x,y),n) \\
\nonumber&  \quad \quad =  u(\eta'(x,y),n+1)\times  u(\eta'(x,y),n-1)\\
\nonumber&  \quad \quad \quad \quad \quad +  u(\eta'(x+\frac{1}{2},y-\frac{1}{2}),n)\times  
 u(\eta'(x-\frac{1}{2},y+\frac{1}{2}),n).
\end{align}
\end{prop}

If $y=x-1$, we get the same kind of character formulas, but for the   $\bar u(\eta'(x,y),n)$ : 

\begin{equation}\label{Zu4}
\bar u(\eta'(x,x-1),n)= (-1)^{n+1} \det((\tilde \eta'_{i,j})_{1\leq i,j \leq n} ),
\end{equation}
where   $\tilde \eta'_{i,j}= \tilde \eta' (x+\frac{n-1}{2}-i+1,y+\frac{n-1}{2}-j+1 )$, and $\tilde \eta'$  
denotes the coherent
family coinciding with $\eta$ when $x-y$ is positive, as in the real case.

Again  from the Lewis Carroll identity, we deduce  the following (with $2n$ in place of $n$):

\begin{align}\label{groscalcul}
& \nu^{-\frac{1}{2}}  \bar u(\eta'(x,x-1) ,2n)\times \nu^{\frac{1}{2}} \bar u(\eta'(x,x-1),2n) \\
\nonumber&  \quad \quad =  \bar u(\eta'(x,x-1),2n+1)\times  \bar u(\eta(x,x-1),2n-1)\\
\nonumber&  \quad \quad \quad \quad \quad +  \bar u(\eta'(x+\frac{1}{2},x-\frac{1}{2}),2n)\times  
 \bar u(\eta(x-\frac{1}{2},x-\frac{3}{2}),2n).
\end{align}

The $\bar u(\eta'(.,.),.)$ can be expressed as products of  $u(\eta'(.,.),.)$, explicitly: 
\begin{align*}
 \bar u(\eta'(x,x-1) ,2n)&=  u(\eta'(x+\frac{1}{2},x-\frac{1}{2}) ,n)\times  u(\eta'(x-\frac{1}{2},x-\frac{3}{2}) ,n)\\
 \bar u(\eta'(x,x-1) ,2n+1)&=  u(\eta'(x,x-1) ,n+1)\times  u(\eta'(x,x-1) ,n)
\end{align*}

Substituting in  this in (\ref{groscalcul}), and using the fact that the ring $\caR$ is a domain, we find that : 

\begin{prop}
\begin{align}\label{boutH1}
& \nu^{-1} u(\eta'(x,x-1) ,n)\times \nu  u(\eta'(x,x-1),n) \\
\nonumber&  \quad \quad =  u(\eta'(x,x-1),n+1)\times  u(\eta'(x,x-1),n-1)\\
\nonumber&  \quad \quad \quad \quad \quad +  u(\eta'(x+\frac{1}{2},x-\frac{1}{2}),n)\times  
 u(\eta'(x-\frac{1}{2},x-\frac{3}{2}),n).
\end{align}
\end{prop}

\bigskip

\section{Compatibility and further comments}

Let $F$ be a local field (archimedean or non archimedean of any characteristic) and $A$ a central division
 algebra of dimension $d^2$ over $F$ (if $F$ is archimedean, then $d\in\{1,2\}$). If $g\in G^F_{nd}$ is a 
regular semisimple element, we say that $g$ {\it transfers} if there exists an element $g'$ of $G_n^A$ 
which corresponds to $g$ (see Section \ref{JL1}). Then $g$ transfers if and only if its characteristic
 polynomial breaks into a product of irreducible polynomials of degrees divisible by $d$. We say that 
$\pi\in \caR(G_{nd}^F)$ is {\it $d$-compatible} if $\JL(\pi)\neq 0$. Otherwise stated, $\pi$ is 
$d$-compatible if and only if its character does not identically vanish on the set of elements of
 $G_{nd}^F$ which transfer. This justify the dependence of the definition only on $d$
 (and not on $D$). We then have the following results:

\begin{prop}\label{xxx1}
Let $\pi_i\in \Irr_{n_i}^F$, $1\leq i\leq k$, with $\sum_i n_i=n$. Then $\pi_1\times\pi_2\times...\times \pi_k$ 
is $d$-compatible if and only if for all $1\leq i\leq k$, $d$ divides $n_i$ and $\pi_i$ is $d$-compatible.
\end{prop}

\pr If an element $g\in G_n^F$ is conjugated with an element of a Levi subgroup of $G_n^F$, 
say $(g_1,g_2,...,g_k)\in G_{(n_1,n_2,...,n_k)}$ with $g_i\in G_{n_i}^F$, then the characteristic
 polynomial of $g$ is the product of the characteristic polynomials of $g_i$. It follows that,
 if $g$ is semisimple regular, it transfers if and only if $d|n_i$ for all $i$ and $g_i$ transfers.

It is a general fact that for a fully induced representation of a group $G$ from a Levi subgroup $M$, 
the character is zero on regular semisimple elements which are not conjugated in $G$ to some element in  $M$.
 Moreover, one has a precise formula of the character of the fully induced representation in terms of the
 character of the inducing representation (see \cite{vD} and \cite{Cl}, Proposition 3 for non-archimedean
 $F$, \cite{Kn}  Section 13,  for archimedean $F$).
The proposition follows.\qed

We define now an order $<<$ finer than the Bruhat order on $<$ on $\Irr^A_n$. 
If $\pi=\Lg(\delta_1,\delta_2,...,\delta_k)$ and $\pi'=\Lg(\delta'_1,\delta'_2,...,\delta'_{k'})$
 are in $\Irr^A_n$, we set $\pi<<\pi'$ if 
$$\Lg(\ccc^{-1}(\delta_1),\ccc^{-1}(\delta_2),...,\ccc^{-1}(\delta_k))
<\Lg(\ccc^{-1}(\delta'_1),\ccc^{-1}(\delta'_2),...,\ccc^{-1}(\delta'_{k'}))$$ in $\Irr_{nd}^F$.

\begin{prop}\label{xxx2}
Let $\delta_i\in D_{n_i}^F$, $1\leq i\leq k$. Assume for all $1\leq i\leq k$ we have $d|n_i$,
 and set $\delta'_i={\bf C}(\delta_i)\in D_{n_i}^A$. Then $\Lg(\delta_1,\delta_2,...,\delta_k)$ is compatible
 and one has:
$$\JL(\Lg(\delta_1,\delta_2,...,\delta_k))=(-1)^{nd-n}\Lg(\delta'_1,\delta'_2,...,\delta'_k)+\sum_{j\in J} m_j\pi'_j$$
where $J$ is empty or finite, $m_j\in \bbZ ^*$, $\pi'_j\in \Irr_{\sum n_i}^A$ and 
$\pi'_j<< \Lg(\delta'_1,\delta'_2,...,\delta'_k)$ for all $j\in J$.
\end{prop}

\pr One applies  Theorem \ref{multmat} and an induction on the number of representations smaller than
  $\Lg(\delta_1,\delta_2,...,\delta_k)$. See \cite{Ba3}, Proposition 3.10.
\qed

\begin{prop}\label{xxx3}
If $\delta\in D_n^F$, set $\deg(\delta)=n$ and let $l(\delta)$ be the length of ${\bf Supp}(\delta)$ 
(notice that $l(\delta)|\deg(\delta)$). Then 

a) $u(\delta,k)$ is $d$-compatible if and only if either $d| \deg(\delta)$ or 
$d| k\frac{\deg(\delta)}{l(\delta)}$. 

b) there exists $k_\delta\in \n^*$ such that $u(\delta,k)$ is $d$-compatible if and only if 
$k_\delta |k$. Moreover, $k_\delta |d$.
\end{prop}

\pr a) is  in Section 3.5 of  \cite{Ba1} for the non-archimedean case. It follows from  Theorem 
\ref{archimtransfer} in the archimedean case. 

b) follows easily from a). For the archimedean (non trivial {\sl i.e.} $A=\bbH$) case, $d=2$ and the transfer 
theorem \ref{archimtransfer} shows that 

--- if $\deg(\delta)=2$, then $u(\delta, k)$ is $2$-compatible for all $k$ (hence $k_\delta=1$) and 

--- if $\deg(\delta)=1$ then 
$u(\delta,k)$ is $2$-compatible if (and {\it only} if, because of the dimension of $G_k^F$) $k$ is even 
(hence $k_\delta=2$).\qed

\medskip

Let $\gamma$ be an irreducible generic unitary representation of $G^F_{n}$. As $\gamma$ is generic, 
it is fully induced from an essentially square integrable representation (\cite{Ze} for non 
 archimedean fields, section \ref{GEN} for archimedean fields).
Then as $\gamma$ is unitary, thanks to the classification of the unitary spectrum (\cite{Ta1}, \cite{Vog2} 
and Section \ref{GEN} the present paper), $\gamma$ is an irreducible product 
$\sigma_1\times \sigma_2\times ...\times\sigma_p\times \pi_1\times \pi_2\times ...\times\pi_l$, 
where, for $1\leq i\leq p$, $\sigma_i\in D^{u,F}$, and, for $1\leq j\leq l$, $\pi_j=\pi(\delta_j,1;\alpha_j)$ 
for some $\delta_j\in D^{u,F}$ and some $\alpha_i\in ]0,\frac{1}{2}[$. 

Using the Langlands classification, it is easy to see that the representation 
$$\nu^{\frac{k-1}{2}}\gamma\times  \nu^{\frac{k-1}{2}-1}\gamma\times ...\times \nu^{-\frac{k-1}{2}}\gamma$$
has a unique quotient $u(\gamma,k)$, and one has
$$u(\gamma,k)=u(\sigma_1,k)\times u(\sigma_2,k)\times ...\times u(\sigma_p,k)\times 
\pi(\delta_1,k;\alpha_1) 
\times \pi(\delta_2,k;\alpha_2)\times ...
 \times \pi(\delta_l,k;\alpha_l)$$
(see for instance \cite{Ba3} Section 4.1). The local components of cuspidal automorphic representations of
 $\GL_n$ over adeles of global fields are unitary generic representations (\cite{Sha}).
According to the classification of the residual spectrum (\cite{MW}), it follows that local component of 
residual automorphic representations of the linear group are of type $u(\gamma,k)$.

\begin{prop}\label{xxx4}
Let $\gamma$ be a unitary generic representation of $G_n^F$ for some $n\in \bbN^*$. There exists 
$k_{\gamma}$ such that $u(\gamma,k)$ is $d$-compatible if and only if $k_\gamma |k$. Moreover, $k_\gamma |d$.
\end{prop}

\pr The (easy) proof  given in \cite{Ba1} Section 3.5 for non-archimedean fields works also 
  for archimedean fields. If $$u(\gamma,k)=u(\sigma_1,k)\times u(\sigma_2,k)\times ...
\times u(\sigma_p,k)\times \pi(\delta_1,k;\alpha_1) 
\times \pi(\delta_2,k;\alpha_2)\times ... \times \pi(\delta_l,k;\alpha_l),$$
then $u(\gamma,k)$ is $d$-compatible if and only if all the $u(\sigma_i,k)$ and  $u(\delta_j,k)$
 are compatible (Proposition \ref{xxx1}).
Then  Prop. \ref{xxx3} implies Prop. \ref{xxx4}. If $F=\bbR$, $k_\gamma=1$ if and only if
 all the $\sigma_i$ and $\delta_j$ are in $D_2$. If not, $k_\gamma=2$.\qed

\bigskip

\def\aa{{\mathbb A}}
\def\bc{\backslash}
\def\lg{L^2(Z(\aa)G(F)\bc G(\aa))}
\def\o{\omega}
\def\lgg{L^2(Z(\aa)G'(F)\bc G'(\aa);\o)}
\def\lra{\leftrightarrow}
\def\ccc{{\bf C}}
\def\jlr{{\bf JL}_r}
\def\JLr{{\bf LJ}_r}

\def\s{{\mathfrak S}}
\numberwithin{equation}{section}
\leftskip -1cm
\rightskip -1cm
\def\a{\alpha}
\def\b{{\mathcal B}}
\def\d{{\mathcal D}}
\def\e{\varepsilon}
\def\f{{\mathcal F}}
\def\i{{\bf i}}
\def\k{\{1,2,...,k\}}
\def\l{{\bf l}}
\def\n{\mathbb N}
\def\r{{\mathbb R}}
\def\s{\sigma}
\def\z{\mathbb Z}

\def\cc{{\mathbb C}}
\def\ccc{{\bf C}}
\def\lra{\leftrightarrow}
\def\ski{{\sum_{i=1}^k}}
\def\ki{{_{i=1}^k}}
\def\rrr{{\mathcal R}}

\def\JL{{\bf LJ}}
\def\sgls{standard Levi subgroup}
\def\sglss{standard Levi subgroups}
\def\bc{\backslash}
\def\cusp{{\mathcal C}}

\section{Notation for the global case}

Let $F$ be a global field {\it of characteristic zero} and $D$ a
central division algebra over $F$ of dimension $d^2$. Let
$n\in\n^*$. Set $A=M_{n}(D)$. For each place $v$ of $F$ let $F_v$ be
the completion of $F$ at $v$ and set $A_v=A\otimes F_v$. For every
place $v$ of $F$, $A_v$ is isomorphic to $M_{r_v}(D_v)$ for some positive integer
$r_v$ and some central division algebra $D_v$ of dimension $d_v^2$
over $F_v$ such that $r_v d_v=nd$. We  fix once  for all an
isomorphism $A_v \simeq M_{r_v}(D_v)$  and identify these two algebras. We say that $M_n(D)$
{\it is split} at a place $v$ if $d_v=1$. The set $V$ of places
where $M_n(D)$ is not split is finite. For each $v$, $d_v$
divides $d$, and moreover $d$ is the smallest common multiple of the
$d_v$ over all the places $v$. 

Let $G'(F)$ be the group $A^\times=\GL_n(D)$. For every finite
place $v$ of $F$, set $G'_v= A_v^\times= \GL_{r_v}(D_v)$. For every finite place
$v$ of $F$, we set $K_v=\GL_{r_v}(O_v)$,
where $O_v$ is the ring of
integers of $D_v$. 
 Let $\aa$ be the ring of ad\`eles of
$F$. We define the group $G'(\aa)$ of ad\`eles of
$G'(F)$ as the restricted product of the $G'_v$ over all $v$, with respect to the
family of open compact subgroups $K_v$, $v$ finite. 

Let $G'_\infty$ be the direct product of $G'_v$ over the set of infinite places of $F$ and and $G'_f$ the 
restricted product of $G'_v$ over finite places, with respect to the open compact subgroups $K_v$. The group 
$G'(\aa)$ decomposes into the direct product
$$G'(\aa)=G'_\infty\times G'_f.$$

Fix maximal compact subgroups $K_v$ at archimedean places $v$ like before, $K_v=\Or(n),\U(n),\Sp(n)$
 according to $G'_v$ being $\GL_n(\bbR)$, $\GL_n(\bbC)$ or $\GL_n(\bbH)$. Let $K_\infty$ (resp. $K_f$) 
be the compact subgroup of $G_\infty$ (resp. of $G'_f$) which is the direct product of $K_v$ over the 
infinite places (resp. finite places) $v$. Let $K$ be $K_\infty\times K_f$ as a (compact) subgroup of 
$G'(\aa)$. Let ${\mathfrak g}_\infty$ be the Lie algebra of $G_\infty$.

An {\it admissible $G'(\aa)$-module} is a linear space $V$ which is both a
 $({\mathfrak g}_\infty,K_\infty)$-module and a $G'_f$ smooth module such that the actions of 
$({\mathfrak g}_\infty,K_\infty)$ and $G'_f$ commute and for all irreducible equivalence class 
of continuous representations $\pi$ of $K$ the $\pi$ isotypic component of $V$ is of finite 
dimension. It is {\it irreducible} if it has no proper sub $G'(\aa)$-module, and {\it unitary}
 if admits a Hermitian product which is invariant under both actions of $({\mathfrak g}_\infty,K_\infty)$ 
and $G'_f$.

If $V$ is an admissible $G'(\aa)$-module, then $V$ is isomorphic with a tensor product 
$V_\infty\otimes V_f$, where $V_\infty$ is an irreducible  $({\mathfrak g}_\infty,K_\infty)$-module 
and $V_f$ is an irreducible smooth representation of $V_f$.

If $(\pi,H)$ is a unitary irreducible admissible $G_f$-module, then $\pi$ breaks into a restricted 
tensor product $\otimes_{v\text{ finite}}\pi_v$ where $\pi_v$ is a unitary irreducible representation 
of $G'_v$ (\cite{JL}, \cite{Langl}, \cite{GGPS} or \cite{Fl}). For almost all $v$, $\pi_v$ has 
a fixed vector under the maximal compact subgroup $K_v$. Such a representation is called {\it spherical}. 
The $\pi_v$ are determined by $\pi$. 
Such a $\pi_v$ is called the {\it local component} of $\pi$ at the place $v$. The set of local 
components $\pi_v$ determines $\pi$.

Let $Z(F)$ be the center of $G'(F)$ and, for every place $v$, let $Z_v$
be the center of $G'_v$. Then we identify the center $Z(\aa)$ of $G'(\aa)$ with the restricted product 
of the $Z_v$, with respect to the open compact subgroups $Z_v\cap K_v$ at finite places. For any finite $v$, 
we fix a Haar measure $dg_v$ on $G'_v$ such that the volume of $K_v$ is one, and a Haar measure $dz_v$ on $Z_v$ 
such that the volume of $Z_v\cap K_v$ is one. The set of measures $\{dg_v\}_{v\ \text{finite}}$ induce a well 
defined Haar measure on the locally compact group  $G'_f$ and $\{dz_v\}_{v\ \text{finite}}$ induce a well
 defined measure on its center (see for instance \cite{RV} where measures on restricted products are explained).

For the archimedean groups we chose
Duflo-Vergne's normalization, defined as follows: let $G$ be a reductive
 group (complex or real), and pick a  $G$-invariant symmetric,
non-degenerate bilinear form $\kappa$ on the Lie algebra $\frg$. Then $\frg$ will be
endowed with the Lebesgue measure $dX$ such that the volume of a
parallelotope supported by a basis $\{X_1, \ldots, X_n\}$ of $\frg$ is
equal to $|\det(\kappa(X_i,X_j))|^{\frac{1}{2}}$ and $G$ will be endowed
with the Haar measure tangent to $dX$. If $G'$ is a closed subgroup of $G$,
such that $\kappa$ is non-degenerate on its Lie algebra $\frg'$, we endow $G'$ with the
Haar measure determined by $\kappa$ as above. This gives measures on $G'_\infty$ and its center. 

We fix now the measure $dg$ on $G'(\aa)=G'_\infty\times G'_f$ (resp. $dz$ on $Z(\aa)$) which is 
the product of measures chosen before for the infinite and the finite part. We fix a measure on
 $Z(\aa)\bc G'(\aa)$ which is the quotient measure $dz\bc dg$.

We see $G'(F)$ as a subgroup of
$G'(\aa)$ via the diagonal embedding. As $G'(F)\cap Z(\aa)\bc
G'(F)$ is a discrete subgroup of $Z(\aa)\bc G'(\aa)$, $dz\bc dg$ defines a measure on the
quotient space $Z(\aa)G'(F)\bc G'(\aa)$. The measure of the space
$Z(\aa)G'(F)\bc G'(\aa)$ is finite. 

Fix a unitary smooth character $\o$ of $Z(\aa)$, trivial on $Z(F)$.

Let $\lgg$ be the space of classes of functions $f$ defined on
$G'(\aa)$ with values in $\cc$ such that

i) $f$ is left invariant under $G'(F)$,

ii) $f$ satisfies $f(zg)=\o(z)f(g)$ for all $z\in Z(\aa)$ and almost
all $g\in G'(\aa)$,

iii) $|f|^2$ is integrable over $Z(\aa)G'(F)\bc G'(\aa)$.\\

Let $R'_\o$ be the representation of $G'(\aa)$ in $\lgg$ by right translations. As explained in
\cite{BJ}, each irreducible subspace of $\lgg$ gives rise to a unique unitary irreducible admissible 
$G'(\aa)$-module. We call such a $G'(\aa)$-module a {\it discrete series of}
$G'(\aa)$.

Every discrete series of $G'(\aa)$ with the central character $\o$
appears in $R'_\o$ with a finite multiplicity (\cite{GGPS}).

Let $R'_{\o,disc}$ be the subrepresentation of $R'_\o$ generated by
the discrete series.  If $\pi$ is a discrete series we call the {\it
multiplicity of $\pi$ in the discrete spectrum} the multiplicity
with which $\pi$ appears in $R'_{\o,disc}$.\\

{\bf Notation.} Fix $n$ and $D$ as before. The same constructions work obviously starting with $A=\GL_{nd}(F)$ instead of 
$A=\GL_n(D)$. We  denote $G(\aa)$ the group of invertible elements of $A$  
and modify all notations accordingly.\\

\section{Second insight of some local results}

We would like to point out that some of the archimedean results described in this paper may be proved by 
global methods and local tricks as in the non-archimedean case (\cite{Ba1} and \cite{Ba3}), avoiding any 
reference to cohomological induction. These are $U(1)$ for $\GL(n,\bbH)$, the fact that products of 
representations in $\caU_\bbH$ are irreducible and the Jacquet-Langlands transfer of unitary
 representations (using $U(0)$ for $\GL(n,\bbR)$ - \cite{Bar} - but not  on $\GL(n,\bbH)$)
We sketch here these proofs.

\subsection{  $U(1)$  and the transfer of $u(\delta,k)$} 

Let $\JL:\caR_{2n}^\bbR\to \caR_{n}^\bbH$ be  the morphism between Grothendieck
 groups extending the classical Jacquet-Langlands correspondence for square integrable 
representations (Section \ref{JL1}). We give here a second proof of the 

\begin{prop}\label{newproof}
\rm{(a)} If $\chi\in D_1$, then $\JL(u(\chi,2n))=\chi'_n$.

\rm{(b)} If $\delta\in D_2$ and $\delta'=\bf{C}(\delta)$, then $\JL(u(\delta,n))=(-1)^n \bar{u}(\delta',n)$.

\rm{(c)} The statement $U(1)$, i.e. {\rm $u(\delta',n)$ are unitary}, is true for $\GL(n,\bbH)$.
\end{prop}

The first assertion (a) is obvious since $u(\chi,2n)=\chi_{2n}$ and the equality of characters may
 be checked directly. To prove (c), recall we have

\begin{equation}\label{new}
\JL(u(\delta,n))=(-1)^n(\bar{u}(\delta',n)+\sum_{i=1}^k a_i u_i),
\end{equation}
where the $u_i$ are irreducible non-equivalent representations of $\GL(n,\bbH)$, non equivalent to
 $\bar{u}(\delta',n)$, and $a_i$ are non-zero integers (Proposition \ref{xxx2}). 

We now claim that all the irreducible representations on the right hand side of the equality are 
unitary and the $a_i$ are all positive. One may proceed like in \cite{Ba3}: choose a global field $F$ 
and a division algebra $D$ over $F$ such that, if $G'(\aa)$ is the adele group of $D^\times$, we have 
$G'_v=GL_n(\bbH)$ for some place $v$. As $\delta\in D_2$, there exists a cuspidal representation $\rho$ 
of $G(\aa)=GL_{2n}(\aa)$ such that $\rho_v=\delta$. According to the classification of the residual 
spectrum for $G(\aa)$ (\cite{MW}) there exists a residual representation $\pi$ of $G(\aa)$ such that 
$\pi_v=u(\delta,n)$.  Comparing then the trace formula from \cite{AC} (or the simple trace formula from
 \cite{Ar}) of $G(\aa)$ and $G'(\aa)$, one gets using standard simplifications and multiplicity one on 
the $G(\aa)$ side a local formula $\JL(u(\delta,n))=\pm\sum_{j=1}^k b_j w_i$ where the $b_j$ are 
multiplicities of representations - hence positive, and $w_j$ are local component of global discrete series 
- hence unitary. By linear independence of characters on $\GL(n,\bbH)$, this formula is the same as the formula
(\ref{new}) which implies in particular $\bar{u}(\delta',n)$ is unitary (see \cite{Ba3}, Cor. 4.8(a)).
 This implies the assertion $U(1)$, since when $\delta'$ is not a character one has
 $\bar{u}(\delta',n)=u(\delta',n)$, while when $\delta'$ is a (unitary) character
 we know $u(\delta',k)$ is the unitary character 
 $\delta'\circ RN$.  So (c) is proved.
 
We now prove (b). We want to prove that on the left hand side of the equality (\ref{new}) there
 is just one term, $\bar{u}(\delta',n)$. If $\pi$ is an irreducible unitary representation of $\GL(n,\bbR)$ we say
 $\pi$ is
 {\it semirigid} if it is a product of representations $u(\delta,k)$. We already showed in the previous
 paragraph that all these representations $u(\delta, k)$ correspond by $\JL$ to zero or a sum of unitary representations.
 As $\JL$ commutes with products and a product of irreducible  unitary representations is a sum of irreducible 
 unitary  representations, it follows that   any sum of semirigid irreducible 
unitary representation of some $\GL(2n,\bbR)$ correspond to zero or a sum of unitary representations 
of $\GL(n,\bbH)$. The relation (\ref{new}) shows now that for all $\alpha\in\r$,
 $\JL(\pi(\delta,n;\alpha))=\nu'^\alpha(\sum_{i=0}^k a_i u_i)\times \nu'^{-\alpha}(\sum_{i=0}^k a_i u_i)$
 where $a_0=1$, $u_0=\bar{u}(\delta',n)$.  When $\alpha=\frac{1}{2}$ on the left hand side of
 the equality we obtain a sum of semirigid unitary representations  (see 
Proposition \ref{boutR} for precise formula), so on the right hand side we should have a sum of 
unitary representations. But this is impossible as soon as the sum $\sum_{i=1}^k a_i u_i$ contains
a representation $u_1$, since then the mixed product $\nu'^{-\frac{1}{2}}u_0
\times \nu'^{\frac{1}{2}}u_1$ contains a non hermitian subquotient (the ``bigger'' one for the Bruhat order for example). 
This shows there is only  one $u_i$, $i=0$, and so $\LJ(u(\delta,n))=(-1)^n\bar{u}(\delta',n)$. 
\qed

\subsection{ Irreducibility  and  transfer of all unitary representations}

We know now that the representations in $\caU_\bbH$ are all unitary. To show that their 
 products remain irreducible, we may use the irreducibility trick in \cite{Ba1}, Proposition 2.13
 which reduces the problem to show that  $u(\delta',k)\times u(\delta',k)$ is irreducible for all discrete series 
$\delta'$ of $\GL(1,\bbH)$  and all $k \in \bbN^*$.   Let $\delta$ be a square
 integrable representation of $\GL(2,\bbR)$ such that $\LJ(\delta)=\delta'$. It follows that we have the
 equality $\LJ(u(\delta,k)\times u(\delta,k))=\bar{u}(\delta',k)\times \bar{u}(\delta',k)$. On the left hand
 side we have the irreducible representation $M=u(\delta,k)\times u(\delta,k)$.  On the right hand
 side we have a sum of unitary representations, the product $M'=\bar{u}(\delta',k)\times \bar{u}(\delta',k)$ 
 (we already know $\bar{u}(\delta',k)$ is unitary), 
 which we want to show has actually a single term. Apply the same $\alpha$ trick like before : we know that
 $\pi(M,\alpha)$ corresponds to $\pi(M',\alpha)$. For $\alpha=\frac{1}{2}$ the first representation
 breaks into a sum of semirigid unitary representations, while the second is a sum containing non unitary
 representations unless $M'$ contains a single term. Notice that the Langlands quotient theorem and 
$U(4)$ guarantee $M'$ has a subquotient which appears with multiplicity one, so either $M'$ is a sum containing
 two different terms, or is irreducible. So the square of $\bar{u}(\delta',k)$ is irreducible for all $k$. If $\delta'$ is not a character, then $u(\delta',k)=\bar{u}(\delta',k)$ so the square of $u(\delta',k)$ is irreducible. If $\delta'$ is a character then we saw $\bar{u}(\delta',2k+1)=u(\delta',k)\times u(\delta',k+1)$ and the result follows again.

This implies now: {\it if $u$ is an irreducible unitary representation of $\GL(2n,\bbR)$, then 
$\LJ(u)$ is either zero, or plus or minus irreducible unitary representation of $\GL(n,\bbH)$}.\\


The proofs here are based on the trace formula and do not involve cohomological induction. 
However, the really difficult result is $U(0)$ on $\GL(n,\bbH)$, and it does.\\

\section{Global results}

\subsection{Global Jacquet-Langlands, multiplicity one and strong multiplicity one for inner forms}

\def\g{{\bf G}}

For all $v\in V$, denote $\JL_v$ (resp. $|\JL|_v$) the
correspondence $\JL$ (resp. $|\JL|$), as defined in Sections \ref{JL1} and  \ref{JaLa}, applied to $G_{v}$ and $G'_{v}$.

If $\pi$ is a discrete series of $G(\aa)$, we say $\pi$ is $D${\it -compatible} if, for all
$v\in V$, $\pi_v$ is $d_v$-compatible. Then $\JL(\pi_v)\neq 0$ and
$|\JL|_v(\pi_v)$ is an irreducible representation of $G'_{n}$.
\def\f{{\bf G}}

Here are the Jacquet-Langlands correspondence and the multiplicity one theorems for $G'(\aa)$
 (already known for $G(\aa)$: \cite{Sha}, \cite{PS}).

\begin{thm}\label{correspondence}
{\rm (a)} There exists a unique map $\f$ from the set of discrete series of $G'(\aa)$ into the set 
of discrete series of $G(\aa)$ such
that $\f(\pi')=\pi$ implies
$|\JL|_v(\pi_v)=\pi'_v$ for all places $v\in V$, and $\pi_v=\pi'_v$
for all places $v\notin V$. The map $\f$ is injective and onto the set of $D$-compatible discrete series of
$G(\aa)$.

{\rm (b)} The multiplicity of every discrete series of $G'(\aa)$ in the discrete spectrum is $1$. If two 
discrete series of $G'(\aa)$ have isomorphic local component at almost every place, then they are equal.
\end{thm}

The proof is the same as the proof of  Theorem 5.1 in \cite{Ba1} with the following minor changes: 
Lemma 5.2 \cite{Ba1} is obviously still true when the inner form is not split at infinite places 
using the Proposition \ref{xxx1} here. For the finiteness property quoted in \cite{Ba1}, p. 417 as [BB],
 one has to replace this reference with \cite{Ba5}, where the case of ramified at infinite places inner 
form is addressed. We do not need here the claim (d) in \cite{Ba1}, Theorem 5.1 which is now a particular
 case of Tadi\'c classification of unitary representation for inner forms. At the bottom of
 pages 417 and 419 in \cite{Ba1}, the independence of characters on a product of connected $p$-adic 
groups is used. Here the product involve also real, sometimes non connected groups like $\GL(n,\bbR)$. 
The linear independence of characters on each of these $\GL_n$  is enough to ensure the linear
 independence of characters on the product, as at infinite places representations are Harish-Chandra
 modules so for all these groups, real or $p$-adic, irreducible representations correspond to 
irreducible modules on a well chosen algebra with idempotents.

As in \cite{Ba1}, the hard core of the proof is the powerful equality 17.8 from \cite{AC}
 (comparison of trace formulae of $G(\aa)$ and $G'(\aa)$).\\
\ \\
\ \\

Let us show now the classification of cuspidal representations of $G'(\aa)$ in terms of cuspidal representations
 of $G(\aa)$. Let $\nu$ (resp. $\nu'$) be the global character of $G(\aa)$ (resp. $G'(\aa)$) given by the product
 of local characters like before (i.e. absolute value of the reduced norm). Recall that, according to
 Moeglin-Waldspurger classification, every discrete series $\pi$ of $G(\aa)$ is the unique irreducible quotient 
of an induced representation $\nu^{\frac{k-1}{2}}\rho \times \nu^{\frac{k-3}{2}}\rho\times... 
\times\nu^{-\frac{k-1}{2}}\rho$ where $\rho$ is cuspidal. Then $k$ and $\rho$ are determined 
by $\pi$, so $\pi$ is cuspidal if and only if $k=1$. We set $\pi=MW(\rho,k)$.

\begin{prop}
\rm{(a)} Let $n\in \bbN^*$ and let $\rho$ be a cuspidal representation of $G_n(\aa)$. Then there exists 
$k_\rho$ such that, if $k\in\n^\times$, then $MW(\rho ,k)$ is $D$-compatible if and only if $k_\rho|k$.
 Moreover, $k_\rho|d$. 

\rm{(b)} Let $\pi'$ be a discrete series of $G'(\aa)$ and $\pi=\g(\pi')$. Then $\pi'$ is cuspidal if
 and only if $\pi$ is of the form $MW(\rho,k_\rho)$. 

\rm{(c)} Let $\rho'$ be a cuspidal representation of some $G'_n(\aa)$. Write $\g(\rho')=MW(\rho,k_\rho)$
 and then set $\nu_{\rho'}=\nu^{k_\rho}$. For every $k\in \bbN^*$, the induced representation 
$$\nu_{\rho'}^{\frac{k-1}{2}}\rho' \times \nu_{\rho'}^{\frac{k-3}{2}}\rho'\times... 
\times\nu_{\rho'}^{-\frac{k-1}{2}}\rho'$$
has a unique irreducible quotient which we will denote $MW'(\rho',k)$. It is a discrete series
 and all discrete series are obtained from some cuspidal $\rho'$ like that. If $\g(\rho')=MW(\rho,k_\rho)$
we have $\g(MW'(\rho',k))=MW(\rho, kk_\rho)$. 
\end{prop}

{\bf Proof.} 
(a) This follows from the Proposition \ref{xxx4} and the fact that for all $v\in V$ $d_v|d$.

(b) This is the proposition 5.5 in \cite{Ba1}, with  ``cuspidal'' in place of ``basic cuspidal''
 thanks to Grbac's appendix. Both the proof of the claim and the proof in the appendix work the same way  here.

(c) When $G'_n(\aa)$ is split at infinite places, this is the claim (a) of Proposition 5.7 in \cite{Ba1}.
 We follow the same idea which reduces the problem to local computation. As \cite{Ba1}
makes use of Zelevinsky  involution, we have to give here a  proof in the archimedean case 
(where the involution doesn't exist). First, to show 
that the induced representation 
$$\nu_{\rho'}^{\frac{k-1}{2}}\rho' \times \nu_{\rho'}^{\frac{k-3}{2}}\rho' \times ... 
\times\nu_{\rho'}^{-\frac{k-1}{2}}\rho'$$
has a constituent which is a discrete series, we will directly show that $\g^{-1}(MW(\rho, kk_\rho))$, 
which is a discrete series indeed, is a constituent of 
$$\nu_{\rho'}^{\frac{k-1}{2}}\rho' \times \nu_{\rho'}^{\frac{k-3}{2}}\rho'\times ...
\times \nu_{\rho'}^{-\frac{k-1}{2}}\rho'.$$
We will show it place by place, local component by local component. Fix a place 
$v$ and let $\gamma$ be the local component of $\rho$ at the place $v$. It is an 
irreducible unitary generic representation, and we know that $u(\gamma,k_\rho)$ transfers.
 Set $\pi=\LJ(u(\gamma,k_\rho))$. What we want to prove is that $\LJ(u(\gamma, kk_\rho))$ 
is a subquotient of $\nu^{k_\rho\frac{k-1}{2}}\pi \times \nu^{k_\rho \frac{k-3}{2}}\pi\times... 
\times\nu^{k_\rho (-\frac{k-1}{2})}\pi$. The unitary generic representation $\gamma$ may be written as
$\gamma=(\times_i\, \sigma_i)\times (\times_j\, \pi (\tau_j,1,\alpha_j))$, with $\sigma_i$ and $\tau_j$
 square integrable representations and $\alpha_j\in ]0,\frac{1}{2}[$. So it is enough to prove
 the result when $\gamma$ is a square integrable representation. Let us suppose $\gamma$ is square 
integrable. To prove that $\pi=\LJ(u(\gamma,k_\rho))$ implies $\LJ(u(\gamma, kk_\rho))$ is a
 quotient of $\nu^{k_\rho\frac{k-1}{2}}\pi \times \nu^{k_\rho \frac{k-3}{2}}\pi\times...
 \times\nu^{k_\rho (-\frac{k-1}{2})}\pi$ we would like to show that the essentially square 
integrable support of the representation $\LJ(u(\gamma, kk_\rho))$ is the union of the square
 integrable support of the representations $\{\nu^{k_\rho(\frac{k-1}{2}-i)}\pi\}_{i\in \{0,1,...,k-1\}}$. 
Then, as the essentially square integrable support of $\times_{i=0}^{k-1}[\nu^{k_\rho(\frac{k-1}{2}-i)}\pi ]$ is in standard
 order,   $\LJ(u(\gamma, kk_\rho))$ will be the unique quotient of the product.
 
If $\gamma$ lives on a group of a size such that it transfers to some $\ccc(\gamma)$, 
then $\pi=\bar{u}(\ccc(\gamma),k_\rho)$, 
$\LJ(u(\gamma, kk_\rho))=\bar{u}(\ccc(\gamma),kk_\rho)$ (\cite{Ba1} Proposition 3.7 (a)
 and second case of transfer in Theorem \ref{archimtransfer} of this paper),
and the result is straightforward. If not, then $u(\gamma, k_\rho)$ verifies the ``twisted'' case
 of transfer \cite{Ba1}, Proposition 3.7 (b) for non archimedean field, first case of Theorem \ref{archimtransfer}
 in this paper for archimedean field. In the non archimedean case, one may compute more
 explicit formulas for the transfer (\cite{Ba1} formula (3.9)) and see that it works. In the archimedean 
case $\gamma$ is a character of $\GL_1(\bbR)$ and so $\pi=\gamma\circ RN_{\frac{k_\rho}{2}}$ 
and $\LJ(u(\gamma, kk_\rho))=\gamma\circ RN_{\frac{kk_\rho}{2}}$.\qed
\ \\
\ \\

Let us recall the uniqueness of the cuspidal support for automorphic representations. According to a result of 
Langlands \cite{LanCor} particularized to our case, we know that any automorphic representation of
 $G'(\aa)$ is a constituent
 of an induced representation of the form $\nu'^{a_1}\rho_1\times \nu'^{a_2}\rho_2\times...
\times \nu'^{a_k}\rho_k$ where $a_i$ are real numbers and $\rho_i$ are cuspidal representations.
 In \cite{JS} the authors prove that, for $G(\aa)$, the couples $(\rho_i,a_i)$ are unique
 which in particular solves the question of the existence of CAP representations. In  \cite{Ba1}
 it is shown the result is true (more or less by transfer) for the more general case $G'(\aa)$, 
if the inner form is split at infinite places. Using the previous  results, the same proof now  works
 with no condition on the infinite places.\\
\ \\

\section{$L$-functions $\epsilon$-factors and transfer}

The fundamental work of Jacquet, Langlands and Godement of $L$-functions and $\epsilon$-factors 
of linear groups on division algebras easily implies the following theorem. What we call $\epsilon'$-factors 
following \cite{GJ} are sometimes called $\gamma$-factors in literature. The value of all functions depend on 
the choice of some additive non trivial
character $\psi$ of $\bbR$ which is not relevant for the results.

\begin{thm}
\rm{(a)} Let $u$ be a $2$-compatible irreducible unitary representation of $\GL_{2n}(\bbR)$ 
and $u'$ the irreducible unitary representation of $\GL_n(\bbH)$ such that $\LJ(u)=\pm u'$. 
Then the $\epsilon'$ factors of $u$ and $u'$ are equal.

\rm{(b)} Let $\delta\in D_2$ and set $\delta'={\bf C}(\delta)$. Then for all $k\in\n^\times$ the 
$L$-functions of $u(\delta,k)$ and $\bar{u}(\delta',k)$ are equal and the $\epsilon$-factors of
 $u(\delta,k)$ and $\bar{u}(\delta',k)$ are equal.

\rm{(c)} If $\chi$ is a character of $\GL(2n,\bbR)$ and $\chi'=\LJ(\chi)$, then the $\epsilon'$-factors 
 of $\chi$ and $\chi'$ are equal.

\end{thm}

{\bf Proof.} If we prove (b) and (c), then (a) follows by the corollary 8.9 from
 \cite{GJ} and classifications of unitary representations in Tadi\'c setting explained in the present paper.

(b) is proved in \cite{JL} for $k=1$. As a particular case of \cite{Jacq} (5.4) page 80, 
the $L$-function (resp. $\epsilon$-factor) of a Langlands quotient $u(\delta,k)$ is the product to the $L$-functions (resp. $\epsilon$-factors)
of 
representations $\nu^{i-\frac{k-1}{2}}\delta$, $0\leq i\leq k-1$. The same proof given there for 
$\GL_{2n}(\bbR)$ works for $\GL_n(\bbH)$ as well, so the case $k=1$ imply the general case. 

(c) In case $\chi$ is the trivial character, this is the corollary 8.10 page 121 in \cite{GJ}. The general case 
follows easily by torsion with $\chi$ (or by reproducing the same proof).\qed



\bibliographystyle{abbrv}

\bibliography{BR3bib}
 
\end{document}